\documentclass[aos,MSNbibl,nameyear,dvips]{arximspdf}
\usepackage{dcolumn}

\doi{10.1214/11-AOS967}
\volume{40}
\issue{1}
\pubyear{2012}
\firstpage{466}
\lastpage{493}

\makeatletter

\newcolumntype{d}[1]{D{.}{.}{#1}}

\renewcommand{\bar}{\overline}

\newcommand{\cal}{\mathcal}

\newcommand{\LLL}{\mathcal{L}}
\newcommand{\field}[1]{\mathbb{#1}}
\newcommand{\R}{\field{R}}
\newcommand{\N}{\field{N}}
\newcommand{\Z}{\field{Z}}
\newcommand{\E}{\field{E}}
\newcommand{\h}[1]{\bolds{#1}}
\newcommand{\hh}[1]{\mathbf{#1}}
\newcommand{\PP}{\mathcal{P}}
\newcommand{\FF}{\mathcal{F}}
\newcommand{\HH}{\mathcal{H}}
\newcommand{\CC}{\mathcal{C}}
\newcommand{\T}{{T}}
\newcommand{\II}{\mathit{II}}
\newcommand{\III}{\mathit{III}}
\newcommand{\dd}{{\mathrm{d}}}

\newtheorem{theorem}{Theorem}
\newtheorem{lemma}[theorem]{Lemma}

\newproclaim{remark}{Remark}

  \let\sv@tabnotetext\tabnotetext
  \let\sv@tabnotemark@fmt\tabnotemark@fmt
   \long\def\legend#1{{\let\tabnote@indent\leavevmode\sv@tabnotetext[]{}{#1}}}

\makeatother

\begin{document}
\begin{frontmatter}

\title{Covariance matrix estimation for stationary time~series\thanksref{T1}}
\runtitle{Estimation of covariance matrices}

\thankstext{T1}{Supported in part by NSF Grants DMS-09-06073 and DMS-11-06970.}

\begin{aug}
\author[A]{\fnms{Han} \snm{Xiao}\corref{}\ead[label=e1]{xiao@galton.uchicago.edu}}
\and
\author[A]{\fnms{Wei Biao} \snm{Wu}\ead[label=e2]{wbwu@galton.uchicago.edu}}
\runauthor{H. Xiao and W. B. Wu}
\affiliation{University of Chicago}
\address[A]{Department of Statistics \\
University of Chicago\\
5734 S. University Ave\\
Chicago, Illinois 60637\\
USA\\
\printead{e1}\\
\hphantom{E-mail: }\printead*{e2}} 
\end{aug}

\received{\smonth{5} \syear{2011}}
\revised{\smonth{12} \syear{2011}}

%
\begin{abstract}
We obtain a~sharp convergence rate for banded covariance matrix
estimates of stationary processes. A~precise order of magnitude is
derived for spectral radius of sample covariance matrices. We also
consider a~thresholded covariance matrix estimator that can better
characterize sparsity if the true covariance matrix is sparse. As our
main tool, we implement Toeplitz [\textit{Math. Ann.} \textbf{70}
(1911) 351--376] idea and relate eigenvalues of covariance matrices to
the spectral densities or Fourier transforms of the covariances. We
develop a~large deviation result for quadratic forms of stationary
processes using \mbox{$m$-dependence} approximation, under the framework of
causal representation and physical dependence measures.
\end{abstract}

%

\end{frontmatter}

\section{Introduction}
\label{secindtroduction}

$\!\!\!$One hundred years ago, in 1911, Toeplitz
obtained a~deep result on eigenvalues of infinite matrices of the form
$\Sigma_\infty= (a_{s-t})_{s, t=-\infty}^\infty$. We say that
$\lambda$ is an eigenvalue of $\Sigma_\infty$ if the matrix
$\Sigma_\infty- \lambda\operatorname{Id}_\infty$ does not have a~bounded
inverse, where $\mathrm{Id}_\infty$ denotes the infinite-dimensional
identity matrix. Toeplitz proved that, interestingly, the set of
eigenvalues is the same as the image set $\{ g(\theta), \theta\in
[0, 2 \pi]\}$, where
%
\begin{equation}
g(\theta) = \sum_{t \in\Z} a_t e^{{i} t \theta} \qquad\mbox{with
} {i}=\sqrt{-1}.
\end{equation}
Note that $g(\theta)$ is the Fourier transform of the sequence
$(a_t)_{t=-\infty}^\infty$. For a~finite $T \times T$ matrix
$\Sigma_{\T} = (a_{s-t})_{1\le s, t \le T}$, its eigenvalues are
approximately equally distributed (in the sense of Hermann Weyl) as
$\{g(\omega_s), s=0, \ldots, T-1 \}$, where $\omega_s = 2 \pi s /
T$ are the Fourier frequencies. See the excellent monograph by
\citet{grenander1958} for a~detailed account.

Covariance matrix is of fundamental importance in many aspects of
statistics including multivariate analysis, principal component
analysis, linear discriminant analysis and graphical modeling.
One can infer dependence structures among variables by
estimating the associated covariance matrices. In the context of
stationary time series analysis, due to stationarity, the
covariance matrix is Toeplitz in that, along the off-diagonals
that are parallel to the main diagonal, the values are constant.
Let $(X_t)_{t \in\Z}$ be a~stationary process with mean $\mu= \E
X_t$, and denote by $\gamma_k =\E[(X_0-\mu)(X_k-\mu)]$, $k \in
\Z$, its autocovariances. Then
%
\begin{equation}\label{eqSgnD28238}
\Sigma_{\T} = (\gamma_{s-t})_{1\leq s,t \leq\T}
\end{equation}
is the autocovariance matrix of $(X_1, \ldots, X_{\T})$. In the rest
of the paper for simplicity we also call~(\ref{eqSgnD28238}) the
covariance matrix of $(X_1, \ldots, X_{\T})$. In time series analysis
it plays a~crucial role in prediction
[\citet{kolmogorov1941}, \citet{wiener1949}], smoothing and
best linear unbiased estimation
(BLUE). For example, in the Wiener--Kolmogorov prediction theory, one
predicts~$X_{\T+1}$ based on past observations $X_{\T}, X_{\T-1},
\ldots.$ If the covariances $\gamma_k$ were known, given observations
$X_1, \ldots, X_{\T}$, the coefficients of the best linear unbiased
predictor $\hat X_{\T+1} = \sum_{s=1}^{\T} a_{\T,s} X_{\T+1-s}$ in
terms of
the mean square error $\| X_{\T+1} - \hat X_{\T+1} \|^2$
are the solution of the discrete Wiener--Hopf equation
\[
\Sigma_{\T} \hh{a}_{\T} = \h{\gamma}_{\T},
\]
where $\hh{a}_{\T} = (a_{\T,1}, a_{\T,2}, \ldots, a_{\T,\T
})^\top
$ and
$\h{\gamma}_{\T} = (\gamma_1, \gamma_2, \ldots, \gamma_{\T
})^\top$,
and we use the superscript $\top$ to denote the transpose of a~vector
or a~matrix. If $\gamma_k$ are not known, we need to estimate them
from the sample $X_1, \ldots, X_{\T}$, and a~good estimate of
$\Sigma_{\T}$ is required. As another example, suppose now $\mu= \E
X_t \not= 0$ and we want to estimate it from $X_1, \ldots, X_{\T}$ by
the form $\hat\mu= \sum_{t=1}^{\T} c_t X_t$, where $c_t$ satisfy the
constraint $\sum_{t=1}^{\T} c_t = 1$. To obtain the BLUE, one
minimizes $\E(\hat\mu- \mu)^2$ subject to $\sum_{t=1}^{\T} c_t = 1$,
ensuring unbiasedness. Note that the usual choice $c_t \equiv1/T$ may
not lead to BLUE. The optimal\vspace*{1pt} coefficients are given by $(c_1, \ldots,
c_{\T})^\top= (\hh{1}^\top\Sigma_{\T}^{-1}
\hh{1})^{-1}\Sigma_{\T}^{-1} \hh{1}$, where
$\hh{1}=(1,\ldots,1)^\top\in\R^\T$; see \citet
{adenstedt1974}. Again a~good estimate of $\Sigma_{\T}^{-1}$ is needed.

Given observations $X_1, X_2, \ldots, X_{\T}$, assuming at the outset
that $\E X_t=0$, we can naturally estimate $\Sigma_{\T}$ via plug-in
by the sample version
%
\begin{equation}\label{eqscovD29945}
\hat\Sigma_{\T} = (\hat\gamma_{s-t})_{1\leq s,t\leq\T}
\qquad\mbox{where }
\hat\gamma_k = {1\over T} \sum_{t=|k|+1}^{\T} X_{t-|k|} X_t.
\end{equation}
To judge the quality of a~matrix estimate, we use the operator
norm. The term ``operator norm'' usually indicates a~class of matrix
norms; in this paper it refers to the $\ell_2/\ell_2$ operator
norm or spectral radius defined as $\lambda({A}) := {\max_{|\hh{x}|=1}}
|A\hh{x}|$ for any symmetric real matrix $A$, where $\hh{x}$
is a~real vector, and\vadjust{\goodbreak} $|\hh{x}|$ denotes its Euclidean norm. For the
estimate $\hat\Sigma_{\T}$ in~(\ref{eqscovD29945}), unfortunately,
because too many parameters or autocovariances are estimated and
the signal-to-noise ratios are too small at large lags,\vspace*{1pt} this
estimate is not consistent. \citet{wuss2009} showed that
$\lambda({\hat\Sigma_{\T} - \Sigma_{\T}}) \not\to0$ in
probability. In
Section~\ref{secinconsistency} we provide a~precise order of
magnitude of $\lambda({\hat\Sigma_{\T} - \Sigma_{\T}})$ and give explicit
upper and lower bounds.

The inconsistency of sample covariance matrices has also been observed
in the context of high-dimensional multivariate analysis, and this
phenomenon is now well understood, thanks to the results from random
matrix theory. See, among others, \citet{marcenko1967},
\citet{bai1993}
and \citet{johnstone2001}. Recently, there is a~surge of interest on
regularized covariance matrix estimation in high-dimensional
statistical inference. We only sample a~few works which are closely related
to our problem. \citet{cai2010}, \citet{bickel2008a} and
\citet{furrer2007} studied the banding and/or tapering methods, while
\citet{bickel2008b} and \citet{karoui2008} considered the
regularization by thresholding. In most of these works, convergence
rates of the estimates were established.

However, none of the techniques used in the aforementioned papers is
applicable in our setting since their estimates require multiple
independent and identically distributed (i.i.d.) copies of random vectors
from the underlying multivariate distribution, though the number of
copies can be far less than the dimension of the vector. In time
series analysis, however, it is typical that only one realization is
available. Hence we shall naturally use the sample
autocovariances. In a~companion paper, \citet{wu2011} established a~systematic theory for ${\cal L}^2$ and ${\cal L}^\infty$ deviations of
sample autocovariances. Based on that, we adopt the regularization
idea and study properties of the banded, tapered and thresholded
estimates of the covariance matrices. \citet{wuss2009} and
\citet{mcmurry2010} applied the banding and tapering methods to the
same problem, but here we shall obtain a~better and optimal
convergence rate. We shall point out that the regularization ideas of
banding and tapering are not novel in time series analysis and they
have been applied in nonparametric spectral density estimation.

In this paper we use the ideas in \citet{toeplitz1911} and
\citet{grenander1958} together with Wu's (\citeyear{wu2005}) recent
theory on
stationary processes to present a~systematic theory for estimates of
covariance matrices of stationary processes. In particular, we shall
exploit the connection between covariance matrices and spectral
density functions and prove a~sharp convergence rate for banded
covariance matrix estimates of stationary processes. Using convergence
properties of periodograms, we derive a~precise order of magnitude for
spectral radius of sample covariance matrices. We also consider a~thresholded
covariance matrix estimator that can better characterize
sparsity if the true covariance matrix is sparse. As a~main technical
tool, we develop a~large deviation type result for quadratic forms of
stationary processes using $m$-dependence approximation, under the
framework of causal representations and physical dependence measures.\vadjust{\goodbreak}

The rest of this article is organized as follows. In Section
\ref{secinconsistency} we introduce the framework of causal
representation and physical dependence measures that are useful for
studying convergence properties of covariance matrix estimates. We
provide in Section~\ref{secinconsistency} upper and lower bounds for
the operator norm of the sample covariance matrices. The convergence
rates of banded/tapered and thresholded sample covariance matrices are
established in Sections~\ref{secbanding} and
\ref{secthresholding}, respectively. We also conduct a~simulation study to compare the finite sample performances of banded
and thresholded estimates in Section~\ref{secsimulation}.
Some useful moment inequalities are collected in Section~\ref{sec6}.
A~large
deviation result about quadratic forms of stationary processes, which
is of independent interest, is given in Section~\ref{secld}. Section
\ref{secconclude} concludes the paper.

We now introduce some notation. For a~random variable $X$ and $p>0$,
we write $X \in\LLL^p$ if $\|X\|_p := (\E|X|^p)^{1/p} < \infty$, and
use $\|X\|$ as a~shorthand for $\|X\|_2$. To express centering of
random variables concisely, we define the operator $\E_0$ as
$\E_0(X):=X-\E X$. Hence $\E_0 (\E_0(X)) = \E_0(X)$. For a~symmetric
real matrix $A$, we use $\lambda_{\min}(A)$ and $\lambda_{\max}(A)$
for its smallest and largest eigenvalues, respectively, and use
$\lambda(A)$ to denote its operator norm or spectral radius. For a~real number $x$, $\lfloor x \rfloor:= \max\{y \in\Z\dvtx y \leq x\}$ denotes
its integer part and $\lceil x \rceil:= \min\{y \in\Z\dvtx y \geq x\}
$. For
two real numbers $x,y$, set $x\vee y =\max\{x,y\}$ and $x\wedge
y:=\min\{x,y\}$. For two sequences of positive numbers $(a_{\T})$ and
$(b_{\T})$, we write $a_{\T} \asymp b_{\T}$ if there exists some
constant $C>1$ such that $C^{-1} \leq a_{\T}/b_{\T} \leq C$ for all
$T$. The letter $C$ denotes a~constant, whose values may vary from
place to place. We sometimes add symbolic subscripts to emphasize
that the value of $C$ depends on the subscripts.\vspace*{-3pt}

\section{Exact order of operator norms of sample covariance
matrices}
\label{secinconsistency}

Suppose $Y$ is a~$p\times n$ random matrix consisting of i.i.d. entries
with mean 0 and variance~1, which could be viewed as a~sample of
size $n$ from some $p$-dimensional population; then $YY^\top/n$ is the
sample covariance matrix. If $\lim_{n \to\infty} p/n=c>0$, then
$YY^\top/n$ is not a~consistent estimate of the population covariance
matrix (which is the identity matrix) in term of the operator
norm. This is a~well-known phenomenon in random matrices literature;
see, for example, \citet{marcenko1967}, Section 5.2 in
\citet{bai2010}, \citet{johnstone2001} and
\citet{karoui2005}. However, the techniques used in those papers are
not applicable here, because we have only one realization and the
dependence within the vector can be quite complicated. Thanks to the
Toeplitz structure of $\Sigma_\T$, our method depends on the
connection between its eigenvalues and the spectral density, defined
by
%
\begin{equation}\label{eqD30222}
f(\theta) = \frac{1}{2\pi}\sum_{k\in\Z} \gamma_k\cos(k\theta).
\end{equation}
The following lemma is a~special case of Section 5.2
[\citet{grenander1958}].\vadjust{\goodbreak}
%
%
\begin{lemma}
\label{thmtoeplitz}
Let $h$ be a~continuous symmetric function on $[-\pi,\pi]$. Denote
by $\underline{h}$ and $\bar{h}$ its minimum and maximum,
respectively. Define $a_k = \int_{-\pi}^{\pi} h(\theta) e^{-{i}
k \theta} \,\dd\theta$ and the $T\times T$ matrix
$\Gamma_{\T}=(a_{s-t})_{1\leq s,t\leq\T}$; then
\[
2\pi\underline{h} \leq\lambda_{\min}(\Gamma_{\T})
\leq\lambda_{\max}(\Gamma_{\T}) \leq2\pi\bar{h}.
\]
\end{lemma}

Lemma~\ref{thmtoeplitz} can be easily proved by noting that
%
\begin{equation}
\label{eqrho}\quad
\hh{x}^\top\Gamma_{\T} \hh{x}
= \int_{-\pi}^\pi| \hh{x}^\top\rho(\theta) |^2 h(\theta) \,\dd
\theta\qquad \mbox{where }
\rho(\theta)=(e^{{i} \theta}, e^{{i}2\theta},
\ldots, e^{{i} \T\theta})^\top.
\end{equation}
%

The sample covariance matrix~(\ref{eqscovD29945}) is closely related
to the periodogram
\[
I_{\T}(\theta) = T^{-1}\Biggl|\sum_{t=1}^{\T} X_t e^{{i}
t\theta}\Biggr|^2 = \sum_{k=1-\T}^{\T-1}
\hat\gamma_k e^{{i} k \theta}.
\]
By Lemma~\ref{thmtoeplitz}, we have $\lambda({\hat\Sigma_{\T}})
\le
\max_{-\pi\le\theta\le\pi} I_{\T}(\theta)$. Asymptotic properties
of periodograms have recently been studied by \citet{peligrad2010}
and \citet{lin2009}. To introduce the result in the latter paper,
we assume that the process $(X_t)$ has the causal representation
%
\begin{equation}
\label{eqcausal}
X_t=g(\varepsilon_t,\varepsilon_{t-1},\ldots),
\end{equation}
where $g$ is a~measurable function such that $X_t$ is well defined,
and $\varepsilon_t, t\in\Z$, are i.i.d. random variables. The framework
(\ref{eqcausal}) is very general [see, e.g.,
\citet{tong1990}] and easy to work with. Let
$\FF^t=(\varepsilon_t,\varepsilon_{t-1},\ldots)$ be the set of innovations
up to time~$t$; we write $X_t=g(\FF^t)$. Let $\varepsilon_t', t\in\Z$,
be an i.i.d. copy of $\varepsilon_t, t\in\Z$. Define $\FF^t_\ast=
(\varepsilon_t, \ldots, \varepsilon_1, \varepsilon_0', \varepsilon_{-1},
\ldots)$,
obtained by replacing~$\varepsilon_0$ in~$\FF^t$ by $\varepsilon_0'$, and
set $X_t'=g(\FF_\ast^t)$. Following \citet{wu2005}, for $p>0$, define
%
\begin{equation}\label{eqtspdm}
\Theta_{p}(m)=\sum_{t=m}^{\infty}\delta_p(t),\qquad m \ge0,\qquad
\mbox{where }
\delta_p(t) = \|X_t-X_t'\|_p.
\end{equation}
In \citet{wu2005}, the quantity $\delta_p(t)$ is called
\textit{physical dependence measure}. We make the convention that
$\delta_p(t) = 0$ for $t<0$. Throughout the article, we assume the
\textit{short-range dependence} condition $\Theta_p := \Theta_p(0) <
\infty$. Under a~mild condition on the tail sum $\Theta_{p}(m)$
(cf. Theorem~\ref{thminconsistency}), \citet{lin2009} obtained the
weak convergence result
%
\begin{equation}
\label{eqlin}
\max_{1\leq s \leq q}
\biggl\{\frac{I_{\T}(2\pi s/T)}{2\pi f(2\pi s/T)}\biggr\}
- \log q \Rightarrow\mathcal{G},
\end{equation}
where $\Rightarrow$ denotes the convergence in distribution,
$\mathcal{G}$ is the Gumbel distribution with the distribution
function $e^{-e^{-x}}$, and $q=\lceil T/2 \rceil-1$. Using this result,
we can provide explicit upper and lower bounds on the operator
norm of the sample covariance matrix.
%
\begin{theorem}
\label{thminconsistency}
Assume $X_t \in\LLL^p$ for some $p>2$ and $\E X_t=0$. If
$\Theta_p(m) = o(1/\log m)$ and
$\min_{\theta}f(\theta)>0$, then
\[
\lim_{\T\to\infty}
P\biggl\{ \frac{\pi[\min_{\theta}f(\theta)]^2 \log
T}{12\Theta_2^2}
\leq\lambda({\hat\Sigma_{\T}}) \leq10\Theta_2^2 \log T \biggr\}
= 1.
\]
\end{theorem}

According to Lemma~\ref{thmtoeplitz}, we know $\lambda_{\max}
(\Sigma_{\T}) \leq2 \pi\max_{\theta}f(\theta)$. As an immediate
consequence of Theorem~\ref{thminconsistency}, there exists a~constant $C>1$ such that
\[
\lim_{\T\to\infty}P[C^{-1}\log T \leq
\lambda({\hat\Sigma_{\T} - \Sigma_{\T}}) \leq C \log T] = 1,
\]
which means the estimate $\hat\Sigma_{\T}$ is not consistent, and the
order of magnitude of $\lambda({\hat\Sigma_{\T} - \Sigma_{\T}})$ is
$\log T$. Earlier,\vspace*{1pt} \citet{wuss2009} also showed that
the plug-in estimate $\hat\Sigma_{\T} = (\hat\gamma_{s-t})_{1\leq s,t
\leq\T}$ is not consistent, namely, $\lambda(\hat\Sigma_{\T} -\break
\Sigma) \not\to 0$ in probability. Our Proposition
\ref{thminconsistency} substantially refines this result by providing
an exact order of magnitude of $\lambda({\hat\Sigma_{\T} - \Sigma})$.

\citet{an1983} showed that under suitable conditions, for linear
processes with i.i.d. innovations,
%
\begin{equation}
\label{eqD306}
\lim_{\T\rightarrow\infty} \max_{\theta} \{I_{\T}(\theta)
/ [2\pi f(\theta) \log T] \}
= 1 \qquad\mbox{almost surely.}
\end{equation}
A~stronger version was found by \citet{turkman1990} for Gaussian
processes. Based on~(\ref{eqD306}), we conjecture that
\[
\lim_{\T\to\infty}
{{\lambda({\hat\Sigma_{\T}})} \over
{2\pi\max_\theta f(\theta) \log T}}
= 1 \qquad\mbox{almost surely.}
\]
\citet{turkman1984} established the following result on the maximum
periodogram of a~sequence of i.i.d. standard normal random variables:
%
\begin{equation}
\label{eqturkman}
\max_{\theta}I_{\T}(\theta) - \log T
- {{\log(\log T)} \over2}
+ {{\log(3/\pi)} \over2} \Rightarrow\mathcal{G}.
\end{equation}
In view of~(\ref{eqlin}) and~(\ref{eqturkman}), we conjecture
that $\lambda({\hat\Sigma_{\T}})$ also converges to the Gumbel
distribution after proper centering and rescaling. Note that the
Gumbel convergence~(\ref{eqturkman}), where the maximum is taken
over the entire interval $\theta\in[-\pi, \pi]$, has a~different
centering term from the one in~(\ref{eqlin}) which is obtained
over Fourier frequencies.

If $Y$ is a~$p \times n$ random matrix consisting of i.i.d. entries,
\citet{geman1980} and \citet{yin1988} proved a~strong convergence
result for the largest eigenvalues of $Y^\top Y$, in the paradigm
where $n,p\rightarrow\infty$ such that $p/n\rightarrow c \in
(0,\infty)$. See also \citet{bai2010} and references therein. If in
addition the entries of $Y$ are i.i.d. complex normal or normal random
variables, \citet{johansson2000} and \citet{johnstone2001} presented
an asymptotic distributional theory and showed that, after proper
normalization, the limiting distribution of the largest eigenvalue
follows the Tracy--Widom law [\citet{tracy1994}]. Again, their methods
depend heavily on the setup that there are i.i.d. copies of a~random
vector with independent entries, and/or the normality assumption, so
they are not applicable here. \citet{bryc2006} studied the limiting
spectral distribution (LSD) of random Toeplitz matrices whose entries
on different sub-diagonals are i.i.d. \citet{solo2010} considered the
LSD of sample covariances matrices generated by a~sample which is
temporally dependent.
\begin{pf*}{Proof of Theorem~\ref{thminconsistency}}
For notational simplicity we let $\underline{f} :=
\min_{\theta}f(\theta)$ and $\bar{f} := \max_{\theta}f(\theta)$. It
follows immediately from~(\ref{eqlin}) that for any $\delta>0$
%
\begin{equation}
\label{eqlowerbound}
\lim_{\T\rightarrow\infty}P\Bigl[\max_{\theta} I_{\T}(\theta)
\geq2(1-\delta)\pi\underline{f} \log T\Bigr] =1.
\end{equation}
The result in~(\ref{eqlin}) is not sufficient to yield an upper
bound of $\max_{\theta}I_{\T}(\theta)$. For this purpose we need to
consider the maxima over a~finer grid and then use
Lemma~\ref{thmzygmund} to extend to the maxima over the whole real
line. Set $j_{\T} = \lfloor T\log T \rfloor$ and $\theta_s := \theta
_{\T,s}
= \pi s/j_{\T}$ for $0 \leq s \leq j_{\T}$. Define $m_{\T} =
\lfloor T^\beta\rfloor$ for some $0<\beta<1$, and $\tilde X_t =
\HH_{t-m_{\T}} X_t = \E(X_t | \varepsilon_{t-m_{\T}},
\varepsilon_{t-m_{\T}+1}, \ldots)$. Let $S_{\T}(\theta) =
\sum_{t=1}^{\T} X_t e^{{i} t\theta}$ be the Fourier transform of
$(X_t)_{1\leq t \leq\T}$, and $\tilde S_{\T}(\theta) =
\sum_{t=1}^{\T} \tilde X_t e^{{i} t\theta}$ for the
$m_{\T}$-dependent sequence $(\tilde X_t)_{1\leq t \leq\T}$. By
Lemma 3.4 of \citet{lin2009}, we have
%
\begin{equation}
\label{eq8}
\max_{0 \leq s \leq j_{\T}} T^{-1/2}|S_{\T}(\theta_s)
- \tilde S_{\T}(\theta_s)| = o_P( (\log T)^{-1/2}).
\end{equation}
Now partition the interval $[1,T]$ into blocks $\mathcal{B}_1,
\mathcal{B}_2, \ldots,\mathcal{B}_{w_{\T}}$ of size $m_{\T}$, where
$w_{\T} = \lfloor T/m_{\T} \rfloor$, and the last block may have size
$m_{\T} \leq|\mathcal{B}|_{w_{\T}} < 2 m_{\T}$. Define the block
sum $U_{\T,k}(\theta) = \sum_{t \in\mathcal{B}_k} \tilde X_t
e^{{i} t \theta}$ for $1 \leq k \leq w_{\T}$. Choose $\beta>0$
small enough such that for some $0 < \gamma< 1/2$, the inequality
%
\begin{equation}
\label{eq9}
1-\beta+\beta p -\gamma(p-1) -1/2<0
\end{equation}
holds. We use truncation and define $\bar U_{\T,k}(\theta) = \E_0
[U_{\T,k}(\theta)\hh{1}_{|U_{\T,k}(\theta)| \leq T^\gamma}
]$. Using similar arguments as equation (5.5) [\citet{lin2009}]
and~(\ref{eq9}), we have
%
\begin{equation}
\label{eq10}
\max_{0\leq s \leq j_{\T}} T^{-1/2}
\Biggl|\sum_{k=1}^{w_{\T}} [U_{\T,k}(\theta_s)
- \bar U_{\T,k}(\theta_s)]\Biggr|
= o_P( (\log T)^{-1/2}).
\end{equation}
Observe that $\bar U_{\T,k_1}(\theta)$ and $\bar U_{\T,k_2}(\theta)$
are independent if $|k_1-k_2|>1$. Let $\mathfrak{R}(z)$ denote the
real part of a~complex number $z$. Split the sum
$\sum_{k=1}^{w_{\T}} \bar U_{\T,k}(\theta)$ into four parts
\[
R_{\T,1}(\theta) = \sum_{k \ \mathrm{odd}}
\mathfrak{R}(\bar U_{\T,k}(\theta)),\qquad
R_{\T,2}(\theta) = \sum_{k \ \mathrm{even}}
\mathfrak{R}(\bar U_{\T,k}(\theta))
\]
and $R_{\T,3}$, $R_{\T,4}$ similarly for the imaginary part of $\bar
U_{\T,k}$. Since $\E|\bar U_{\T,k}(\theta)|^2
\leq\E| U_{\T,k}(\theta)|^2 \leq|\mathcal{B}_k|
\Theta_2^2$, by Bernstein's inequality [cf. \citet{freedman1975}],
\[
\sup_{\theta}P\biggl[|R_{\T,l}(\theta)|
\geq\frac{3\Theta_2}{2\sqrt{2}}\sqrt{T\log T}\biggr]
\leq2\exp\biggl\{-\frac{({9}/{8})\log T}
{1+ 3\Theta_2^{-1}\sqrt{2\log T} T^{\gamma-1/2}}\biggr\}
\]
for $1\leq l \leq4$. It follows that
%
\begin{equation}
\label{eq11}
\lim_{\T\rightarrow\infty} P\Biggl[\max_{0\leq s \leq j_{\T}}
\Biggl|\sum_{k=1}^{w_{\T}}\bar U_{\T,k}(\theta_s)\Biggr|
\geq3\Theta_2 \sqrt{T\log T}\Biggr] = 0.
\end{equation}
Combining~(\ref{eq8}),~(\ref{eq10}) and~(\ref{eq11}), we have
%
\begin{equation}
\label{eq13}
\lim_{\T\rightarrow\infty}
P\Bigl[\max_{0\leq s \leq j_{\T}} I_{\T}(\theta_s)
\leq9.5 \Theta_2^2 \log T\Bigr] =1,
\end{equation}
which together with Lemma~\ref{thmzygmund} implies that
%
\begin{equation}
\label{equpperbound}
\lim_{\T\rightarrow\infty}P\Bigl[\max_{\theta} I_{\T}(\theta)
\leq10 \Theta_2^2 \log T\Bigr] =1.
\end{equation}
The upper bound in Theorem~\ref{thminconsistency} is an
immediate consequence in view of Lemma~\ref{thmtoeplitz}.
For the lower bound, we use the inequality
%
\[
\lambda({\hat\Sigma_{\T}}) \geq\max_{\theta} \{T^{-1}\rho
(\theta)^{\ast}\Sigma_{\T}\rho(\theta)\},
\]
where $\rho(\theta)$ is defined in~(\ref{eqrho}), and
$\rho(\theta)^\ast$ is its Hermitian transpose. Note that
\begin{eqnarray*}
\rho(\theta)^{\ast}\Sigma_{\T}\rho(\theta)
&=& \sum_{s,t=1}^{\T} \hat\gamma_{s-t}
e^{{i} s \theta}e^{-{i} t \theta}\\
&=& \frac{1}{2\pi} \int_{-\pi}^{\pi}
\sum_{s,t=1}^{\T} I_{\T}(\omega)
e^{-{i} (s-t)\omega} e^{{i} (s-t) \theta}
\,\dd\omega\\
& = &\frac{1}{2\pi} \int_{-\pi}^{\pi} I_{\T}(\omega)
\Biggl|\sum_{t=1}^{\T} e^{{i} t(\omega-\theta)}\Biggr|^2
\,\dd\omega\\
&=& \frac{1}{2\pi} \int_{-\pi-\theta}^{\pi-\theta} I_{\T}(\omega
+\theta)
\Biggl|\sum_{t=1}^{\T} e^{{i} t\omega}\Biggr|^2 \,\dd\omega.
\end{eqnarray*}
By Bernstein's inequality on the derivative of trigonometric
polynomials [cf. \citet{zygmund2002}, Theorem 3.13, Chapter X],
we have
\[
\max_{\theta}|I_{\T}'(\theta)| \leq T \cdot\max_{\theta} I_{\T
}(\theta).
\]
Let $\theta_0=\arg\max_{\theta}I_{\T}(\theta)$. Set
$c=(1-\delta)\pi\underline{f} /(10\Theta_2^2)$. By
Lemma~\ref{thmtoeplitz} and~(\ref{eqfact4}), we know $2\pi\bar{f}
\leq\Theta_2^2$, and hence $c \leq1/20$. If $I_{\T}(\theta_0) \geq
2(1-\delta)\pi\underline{f} \log T$ and $ \max_{\theta}
I_{\T}(\theta) \leq10\Theta_2^2 \log T$, then for $|w| \leq c/T$,
we have
\[
I_{\T}(\theta_0+\omega) \geq[2(1-\delta)\pi\underline{f}
-10 c
\Theta_2^2]\log T = (1-\delta)\pi\underline{f} \log T.\vadjust{\goodbreak}
\]
Since $|{\sum_{j=1}^{\T} e^{{i} j\omega}}|^2 \geq10
T^2 / 11$ when $|w| \leq c/T$, it follows that
\begin{eqnarray*}
\rho(\theta_0)^{\ast}\Sigma_{\T}\rho(\theta_0)
&\geq&\frac{1}{2\pi} \cdot(1-\delta)\pi\underline{f} \log T
\cdot\frac{10 T^2}{11} \cdot\frac{2c}{T}\\
&=& \frac{\pi(1-\delta)^2 \underline{f}^2 T\log T}{11\Theta_2^2},
\end{eqnarray*}
which implies that $ \lambda({\hat\Sigma_{\T}}) \geq\pi(1-\delta)^2
\underline{f}^2 \log T/(11\Theta_2^2)$. The proof is completed by
selecting $\delta$ small enough.
\end{pf*}
%
%
\begin{remark}
In the proof, as well as many other places in this article, we
often need to partition an integer interval $[s,t]\subset\N$ into
consecutive blocks $\mathcal{B}_1,\ldots,\mathcal{B}_b$ with the
same size $m$. Since $s-t+1$ may not be a~multiple of $m$, we make
the convention that the last block $\mathcal{B}_b$ has the size $m
\leq|\mathcal{B}_b| <2m$, and all the other ones have the same
size $m$.
\end{remark}

\section{Banded covariance matrix estimates}
\label{secbanding}

In view of Lemma~\ref{thmtoeplitz}, the inconsistency of
$\hat\Sigma_{\T}$ is due to the fact that the periodogram
$I_{\T}(\theta)$ is not a~consistent estimate of the spectral density
$f(\theta)$. To estimate the spectral density consistently, it is very
common to use smoothing. In particular, consider the lag window
estimate
%
\begin{equation}
\label{eqlagwindow}
\hat f_{\T,B_{\T}}(\theta)
= \frac{1}{2\pi} \sum_{k=-B_{\T}}^{B_{\T}} K(k/B_{\T})
\hat\gamma_k \cos(k\theta),
\end{equation}
where $B_{\T}$ is the bandwidth satisfying natural conditions
$B_{\T}\rightarrow\infty$ and $B_{\T}/T\rightarrow0$, and $K(\cdot
)$ is a~symmetric kernel function satisfying
\[
K(0)=1,\qquad |K(x)| \leq1\quad \mbox{and}\quad
K(x)=0 \qquad\mbox{for } |x|>1.
\]
Correspondingly, we define the tapered covariance matrix estimate
\[
\hat\Sigma_{\T,B_{\T}} = \bigl[K\bigl((s-t)/B_{\T}\bigr)
\hat\gamma_{s-t}\bigr]_{1\leq s,t\leq\T}
= \hat\Sigma_{\T} \star W_{\T},
\]
where $\star$ is the Hadamard (or Schur) product, which is formed by
element-wise multiplication of matrices. The term ``tapered'' is
consistent with the terminology of the high-dimensional covariance
regularization literature. However, the reader should not confuse this with
the notion ``data taper'' that is commonly used in time series
analysis. Our tapered estimate parallels a~lag-window estimate of the
spectral density with a~tapered window. As a~special case, if $K(x) =
\hh{1}_{\{|x| \le1\}}$ is the rectangular kernel, then $\hat
\Sigma_{\T,B_{\T}} = ( \hat\gamma_{s-t} \hh{1}_{ \{|s-t| \le
B_{\T
}\}}
)$ is the banded sample covariance matrix. However, this estimate may
not be nonnegative definite. To obtain a~nonnegative definite
estimate, by the Schur product theorem in matrix theory
[\citet{horn1990}], since $\hat\Sigma_{\T}$ is nonnegative definite,
their Schur product $\hat\Sigma_{\T,B_{\T}}$ is also nonnegative
definite if $W_{\T} = [ K((s-t)/B_{\T}) ]_{1 \le s,t
\leq
T}$ is nonnegative definite. The Bartlett or triangular window
$K_B(u) = \max(0, 1-|u|)$ leads to a~positive definite weight matrix
$W_{\T}$ in view of
%
\begin{equation}
\label{eqpos-def-mat}
K_{B}(u) = \int_\R w(x) w(x+u) \,\dd x,
\end{equation}
where $w(x)=\hh{1}_{\{|x|\leq1/2\}}$ is the rectangular window.
Any kernel
function having form~(\ref{eqpos-def-mat}) must be positive
definite. 

Here we shall show that $\hat\Sigma_{\T,B_{\T}}$ is a~consistent
estimate of $\Sigma_{\T}$ and establish a~convergence rate of
$\lambda({\hat\Sigma_{\T,B_{\T}} - \Sigma})$. We first consider the
bias. By the Ger\v{s}gorin theorem\vspace*{1pt} [cf.
\citet{horn1990}, Theorem 6.1.1], we have
$\lambda({\E\hat\Sigma_{\T,B_{\T}} - \Sigma}) \leq b_{\T }$, where
%
\begin{equation}
\label{eqbias}
b_{\T} = 2\sum_{k=1}^{B_{\T}}
\biggl[1-K\biggl(\frac{k}{B_{\T}}\biggr)\biggr]|\gamma_k|
+ \frac{2}{T}\sum_{k=1}^{B_{\T}} k|\gamma_k|
+ 2\sum_{k=B_{\T}+1}^{\T-1}|\gamma_k|.
\end{equation}
%
The first term on the right-hand side in~(\ref{eqbias}) is due to the
choice of the kernel function, whose order of magnitude is determined
by the smoothness of $K(\cdot)$ at zero. In particular, this term
vanishes if $K(\cdot)$ is the rectangular kernel. If $1 - K(u) =
O(u^2)$ at $u = 0$ and $\gamma_k = O(k^{-\beta})$, $\beta> 1$, then
$b_{\T} = O(B_{\T}^{1-\beta})$ if $1 < \beta< 2$, $b_{\T} =
O(B_{\T}^{1-\beta} + T^{-1})$ if $2 < \beta< 3$ and $b_{\T} =
O(B_{\T}^{-2} + T^{-1})$ if $\beta> 3$. Generally, if
$\sum_{k=1}^\infty|\gamma_k| < \infty$, then $b_{\T} \to0$ as
$B_{\T} \to\infty$ and $B_{\T} < T$.

It is more challenging to deal with $\lambda({\hat\Sigma_{\T,B_{\T
}} - \E\hat\Sigma_{\T,B_{\T}}})$. If $X_t \in\LLL^p$ for some $2
< p
\le
4$ and $\E X_t=0$, \citet{wuss2009} obtained 
%
\begin{equation}
\label{eqwu2009}
\lambda({\hat\Sigma_{\T,B_{\T}} - \E\hat\Sigma_{\T,B_{\T}}})
= O_P\biggl( \frac{ B_{\T} \Theta_p^2}{T^{1-2/p} } \biggr).
\end{equation}
The key step of their method is to use the inequality
\[
\lambda({\hat\Sigma_{\T,B_{\T}} - \E\hat\Sigma_{\T,B_{\T}}})
\leq2\sum_{k=0}^{B_{\T}}|K(k/B_{\T})|
|\hat\gamma_k - \E\hat\gamma_k|,
\]
which is also obtained by the Ger\v{s}gorin theorem. 
It turns out that the above rate can be improved by exploiting the
Toeplitz structure of the autocovariance matrix. By
Lemma~\ref{thmtoeplitz},
%
\begin{equation}
\label{eqtoeplitz1}
\lambda({\hat\Sigma_{\T,B_{\T}} - \E\hat\Sigma_{\T,B_{\T}}})
\leq{2\pi\max_{\theta}}|\hat f_{\T,B_{\T}}(\theta)
- \E\hat f_{\T,B_{\T}}(\theta) |.
\end{equation}
%
Since $\hat f_{\T,B_{\T}}(\theta)$ is a~trigonometric polynomial of order
$B_{\T}$, we can bound its maximum by the maximum over a~fine grid. The
following lemma is adapted from \citet{zygmund2002}, Theorem 7.28,
Chapter X.
%
%
\begin{lemma}
\label{thmzygmund}
Let $S(x) = \frac{1}{2}a_0 + \sum_{k=1}^{n} [a_k\cos(kx) + b_k
\sin(kx)]$ be a~trigonometric polynomial of order $n$. For
any $x^\ast\in\R$, $\delta>0$ and $l \geq2 (1+\delta) n$, let
$x_j = x^\ast+2\pi j/l$ for $0 \leq j \leq l$; then
\[
\max_x |S(x)| \leq(1+\delta^{-1})
\max_{0 \leq j \leq l} |S(x_j)|.
\]
\end{lemma}

For $\delta>0$, let $\theta_j = \pi j/(\lceil(1+\delta)B_{\T}
\rceil)$
for $0
\leq j \leq\lceil(1+\delta)B_{\T} \rceil$; then by Lem\-ma~\ref{thmzygmund},
%
\begin{equation}
\label{eqgrid}
{\max_{\theta} }|\hat f_{\T,B_{\T}}(\theta)
- \E\hat f_{\T,B_{\T}}(\theta) |
\leq{(1+\delta^{-1}) \max_{j}}
|\hat f_{\T,B_{\T}}(\theta_j)
- \E\hat f_{\T,B_{\T}}(\theta_j) | .\hspace*{-22pt}
\end{equation}

%
\begin{theorem}
\label{thmbanded}
Assume $X_t \in\LLL^p$ with some $p>4$, $\E X_t=0$, and
$\Theta_p(m)=O(m^{-\alpha})$. Choose the banding parameter $B_{\T}$ to
satisfy $B_{\T} \rightarrow\infty$, and $B_{\T} = O(T^\gamma)$,
for some
%
\begin{equation}
\label{eqdecayrate}
0 < \gamma<1,\qquad \gamma< \alpha p /2\quad \mbox{and}\quad (1-2\alpha
)\gamma< (p-4)/p.
\end{equation}
Then for $b_{\T}$ defined in~(\ref{eqbias}), and $c_p =
(p+4)e^{p/4}\Theta_4^2$,
%
\begin{equation}
\label{eqtapered}
\lim_{\T\rightarrow\infty}
P\Biggl[\lambda({\hat\Sigma_{\T,B_{\T}} - \Sigma_{\T}})
\leq12 c_p \sqrt{\frac{B_{\T}\log B_{\T}}{T}} + b_{\T}\Biggr] =1.
\end{equation}
In particular, if $K(x)=\hh{1}_{\{|x|\leq1\}}$ and $B_{\T}
\asymp
(T/\log T)^{1/(2\alpha+1)}$, then
%
\begin{equation}
\label{eqbanded}
\lambda({\hat\Sigma_{\T,B_{\T}} - \Sigma_{\T}}) =
O_P\biggl[{\biggl(\frac{\log T}{T}
\biggr)^{{\alpha}/({2\alpha+1})}} \biggr].
\end{equation}
\end{theorem}
\begin{pf}
In view of~(\ref{eqbias}), to prove~(\ref{eqtapered}) we only
need to show that
%
\begin{equation}
\label{eqrandom}
\lim_{\T\to\infty} P\Biggl[\lambda({\hat\Sigma_{\T,B_{\T}} - \E
\hat\Sigma_{\T,B_{\T}}})
\leq12 c_p\sqrt{\frac{B_{\T}\log B_{\T}}{T}}
\Biggr] =1.
\end{equation}
By~(\ref{eqtoeplitz1}) and~(\ref{eqgrid}) where we take
$\delta=1$, the problem is reduced to
%
\begin{equation}
\label{eq7}
\lim_{\T\rightarrow\infty}
P\Biggl[ (2\pi) \cdot\max_{j} |\hat f_{\T,B_{\T}}(\theta_j)
- \E\hat f_{\T,B_{\T}}(\theta_j) |
\leq6 c_p\sqrt{\frac{B_{\T}\log B_{\T}}{T}}\Biggr] =1.\hspace*{-15pt}
\end{equation}
By Theorem~\ref{thmldquadratic} (where we take $M=2$), for any
$\gamma<\beta<1$, there exists a~constant $C_{p,\beta}$ such that
%
\begin{eqnarray}
\label{eq32}
&& \max_j P \Biggl[(2\pi) \cdot|\hat f_{\T,B_{\T}}(\theta_j)
- \E\hat f_{\T,B_{\T}}(\theta_j) |
\geq6 c_p \sqrt{\frac{B_{\T}\log B_{\T}}{T}}\Biggr]\nonumber \\
&&\qquad\leq C_{p,\beta} (TB_{\T})^{-p/4} (\log T)
[(TB_{\T})^{p/4}T^{-\alpha\beta p/2}
+ TB_{\T}^{p/2-1-\alpha\beta p/2} + T]\\
&&\qquad\quad{}+ C_{p,\beta} B_{\T}^{-2}.
\nonumber
\end{eqnarray}
If~(\ref{eqdecayrate}) holds, there exist a~$0<\beta<1$ such that
${\gamma-\alpha\beta p/2}<0$ and $(p/4-\alpha\beta
p/2)\gamma-(p/4-1)<0$. It follows that by~(\ref{eq32}),
\begin{eqnarray*}
&& P\Biggl[ {\max_{j}} |\hat f_{\T,B_{\T}}(\theta_j)
- \E\hat f_{\T,B_{\T}}(\theta_j) |
\geq6 c_p\sqrt{\frac{B_{\T}\log B_{\T}}{T}}\Biggr] \\
&&\qquad\leq C_{p,\beta} (\log T)
\bigl[T^{\gamma-\alpha\beta p/2} + T^{1-p/4}
+ T^{(p/4-\alpha\beta p/2)\gamma-(p/4-1)} \bigr]
+ C_{p,\beta}B_{\T}^{-1} \\
&&\qquad= o(1).
\end{eqnarray*}
Therefore,~(\ref{eq7}) holds and the proof of~(\ref{eqtapered})
is complete. The last statement~(\ref{eqbanded}) is an immediate
consequence. Details are omitted.
\end{pf}
%
%
\begin{remark}
\label{rkcentering}
In practice, $\E X_1$ is usually unknown, and we estimate it by
$\bar X_{\T}=T^{-1}\sum_{t=1}^{\T}X_t$. Let
$\hat\gamma_k^{c}=T^{-1}\sum_{t=k+1}^{\T}(X_{t-k}-\bar X_{\T
})(X_t-\bar
X_{\T})$, and $\Sigma_{\T,B_{\T}}^c$ be defined as $\hat\Sigma
_{\T,B_{\T}}$ by
replacing $\hat\gamma_k$ therein by $\hat\gamma_k^c$. Since $\bar
X_{\T} - \E X_1 = O_P(T^{-1/2})$, it is easily seen that
$\lambda(\hat\Sigma_{\T,B_{\T}}-\hat\Sigma_{\T,B_{\T
}}^c)=O_P(B_{\T}/T)$. Therefore,
the results of Theorem~\ref{thmbanded} hold for
$\hat\Sigma_{\T,B_{\T}}^c$ as well.
\end{remark}
%
%
\begin{remark}
In the proof of Theorem~\ref{thmbanded}, we have shown that, as an
intermediate step from~(\ref{eq7}) to~(\ref{eqrandom}),
%
\begin{equation}
\label{eqsdeuniform}
\lim_{\T\to\infty} P \Bigl[{\max_{0 \le\theta\le2 \pi}}
|\hat f_{\T,B_{\T}}(\theta) - \E\hat f_{\T,B_{\T}}(\theta)|
\le6 \pi^{-1} c_p
\sqrt{T^{-1} B_{\T}\log B_{\T} } \Bigr] = 1.\hspace*{-22pt}
\end{equation}
The above uniform convergence result is very useful in spectral
analysis of time series. \citet{shao2007} obtained the weaker
version
\[
\max_{0 \le\theta\le2 \pi} |\hat f_{\T,B_{\T}}(\theta) - \E
\hat
f_{\T,B_{\T}}(\theta) | = O_P \bigl(\sqrt{T^{-1} B_{\T}\log B_{\T
} }\bigr)
\]
under a~stronger assumption that $\Theta_p(m) = O(\rho^m)$ for some
$0<\rho<1$.
\end{remark}
%
%
\begin{remark}
For linear processes, \citet{woodroofe1967} derived the asymptotic
distribution of the maximum deviations of spectral density
estimates. \citet{liuwu2010} generalized their result to nonlinear
processes and showed that the limiting distribution of
\[
\max_{0 \leq j \leq B_{\T}} \sqrt{\frac{T}{B_{\T}}}
\frac{|\hat f_{\T,B_{\T}}(\pi j /B_{\T})
- \E\hat f_{\T,B_{\T}}(\pi j /B_{\T})|}{f(\pi j/B_{\T})}
\]
is Gumbel after suitable centering and rescaling, under stronger
conditions than~(\ref{eqdecayrate}). With their result, and using
similar arguments as Theorem~\ref{thminconsistency}, we can show
that for some constant $C_p$,
\[
\lim_{\T\to\infty}
P\Biggl[C_p^{-1}\sqrt{\frac{B_{\T}\log B_{\T}}{T}}
\leq\lambda({\hat\Sigma_{\T,B_{\T}}-\E\hat\Sigma_{\T,B_{\T}}})
\leq C_p\sqrt{\frac{B_{\T}\log B_{\T}}{T}}\Biggr]=1,
\]
which means that the convergence rate we have obtained in
(\ref{eqrandom}) is optimal.
\end{remark}
%
%
\begin{remark}\label{rem5}
The convergence rate $\sqrt{T^{-1} B_{\T}\log B_{\T}} + b_{\T}$ in
Theorem~\ref{thmbanded} is optimal. Consider a~process $(X_t)$
which satisfies $\gamma_0 = 3$ and when $k>0$,
\[
\gamma_k = \cases{
A^{-\alpha j}, &\quad if $k=A^j$ for some $j\in\N$,\cr
0, &\quad otherwise,}
\]
where $\alpha> 0$ and $A>0$ is an even integer such that
$A^{-\alpha} \le1/5$. Consider the banded estimate
$\hat\Sigma_{\T,B_{\T}}$ with the rectangular kernel. As shown in
the supplementary article [\citet{wu2011s}], there exists a~constant $C>0$ such that
%
\begin{equation}
\label{eqlbd}
\lim_{\T\rightarrow\infty}
P\Biggl[\lambda({\hat\Sigma_{\T,B_{\T}} - \Sigma_{\T}})
\geq C\sqrt{\frac{B_{\T}\log B_{\T}}{T}} + b_{\T}/5\Biggr] =1,
\end{equation}
suggesting that the convergence rate given in~(\ref{eqtapered}) of
Theorem~\ref{thmbanded} is optimal. This optimality property can
have many applications. For example, it can allow one to derive
convergence rates for estimates of $\hh{a}_{\T}$ in the Wiener--Hopf
equation, and the optimal weights $\hh{c}_{\T} = (c_1, \ldots,
c_{\T})^\top$ in the best linear unbiased estimation problem
mentioned in the \hyperref[secindtroduction]{Introduction}.

\end{remark}
%
%
\begin{remark}
We now compare~(\ref{eqwu2009}) and our result. For $p = 4$,
(\ref{eqwu2009}) gives the order $\lambda({\hat\Sigma_{\T,B_{\T}}
- \E\hat\Sigma_{\T,B_{\T}}}) =
O_P(B_{\T}/\sqrt{T})$. Our result~(\ref{eqrandom}) is
sharper by moving the bandwidth $B_{\T}$ inside the square root. We
pay the costs of a~logarithmic factor, a~higher order moment
condition ($p>4$), as well as conditions on the decay rate of tail
sum of physical dependence measures~(\ref{eqdecayrate}). Note that
when $\alpha\ge1/2$, the last two conditions of
(\ref{eqdecayrate}) hold automatically, so we merely need
$0<\gamma<1$, allowing a~very wide range of $B_{\T}$. In
comparison, for~(\ref{eqwu2009}) to be useful, one requires $B_{\T}
= o(T^{1-2/p})$.
\end{remark}
%
%
\begin{remark}
The convergence rate~(\ref{eqwu2009}) of \citet{wuss2009} parallels
the result of \citet{bickel2008a} in the context of high-dimensional
multivariate analysis, which was improved in \citet{cai2010} by
constructing a~class of tapered estimates. Our result parallels the
optimal minimax rate derived in \citet{cai2010}, though the settings
are different.
\end{remark}
%
%
\begin{remark}
Theorem~\ref{thmbanded} uses the operator norm. For the Frobenius
norm see \citet{wu2011} where a~central limit theory for
$\sum_{k=1}^{B_{\T}} \hat\gamma_k^2$ and $\sum_{k=1}^{B_{\T}}
(\hat
\gamma_k - \E\hat\gamma_k)^2$ is established.
\end{remark}

\section{Thresholded covariance matrix estimators}
\label{secthresholding}

In the context of time series, the observations have an intrinsic
temporal order and we expect that observations are weakly
dependent if they are far apart, so banding seems to be natural.
However, if there are many zeros or very weak correlations within
the band, the banding method does not automatically generate a~sparse estimate.\vadjust{\goodbreak} 

The rationale behind the banding operation is sparsity, namely
autocovariances with large lags are small, so it is reasonable to
estimate them as zero. Applying the same idea to the sample covariance
matrix, we can obtain an estimate of $\Sigma_{\T}$ by replacing small
entries in $\hat\Sigma_{\T}$ with zero. This regularization approach,
termed \textit{hard thresholding}, was originally developed in
nonparametric function estimation. It has recently been applied by
\citet{bickel2008b} and \citet{karoui2008} as a~method of covariance
regularization in the context of high-dimensional multivariate
analysis. Since they do not assume any order of the observations,
their sparsity assumptions are permutation-invariant. Unlike their
setup, we still use $\Theta_{p}(m)$ [cf.~(\ref{eqtspdm})] and
%
\begin{equation}
\label{eqstable}\quad
\Psi_{p}(m)=\Biggl(\sum_{t=m}^{\infty}
\delta_{p}(t)^{p'}\Biggr)^{1/p'},\qquad
\Delta_p(m)=\sum_{t=0}^{\infty}
\min\{\CC_p\Psi_{p}(m), \delta_p(t)\}
\end{equation}
as our weak dependence conditions, where $p'=\min(2,p)$ and $\CC_p$ is
given in~(\ref{eqmarzyg}). This is natural for time series analysis.

Let $A_{\T}=2c_p'\sqrt{\log T/T}$, where $c_p'$ is the constant given
in Lemma~\ref{thmmaximumdeviation}. The thresholded sample
autocovariance matrix is defined as
\[
\hat\Gamma_{\T,A_{\T}} = \bigl(\hat\gamma_{s-t}
\hh{1}_{|\hat\gamma_{s-t}|\geq A_{\T}}
\bigr)_{1\leq s,t\leq T}
\]
with the convention that the diagonal elements are never
thresholded. We emphasize that the thresholded estimate may not be
positive definite. The following result says that the thresholded
estimate is also consistent in terms of the operator norm, and
provides a~convergence rate which parallels the banding approach
in Section~\ref{secbanding}. In the proof we compare the
thresholded estimate $\Gamma_{\T,A_{\T}}$ with the banded one
$\Sigma_{\T,B_{\T}}$ for some suitably chosen $B_{\T}$. This is
merely a~way to simplify the arguments. The same results can be proved
without referring to the banded estimates.
%
%
\begin{theorem}
\label{thmthresholding}
Assume $X_t \in\LLL^p$ with some $p>4$, $\E X_t=0$, and
$\Theta_p(m) = O(m^{-\alpha})$, $\Delta_p(m) = O(m^{-\alpha'})$ for
some $\alpha\geq\alpha' > 0$. If
%
\begin{equation}
\label{eqdecayrate2}
\alpha>1/2 \quad\mbox{or}\quad \alpha'p>2,
\end{equation}
then
\[
\lambda({\hat\Gamma_{\T,A_{\T}} - \Sigma_{\T}})
= O_P\biggl[ \biggl(\frac{\log T}{T}
\biggr)^{{\alpha}/({2(1+\alpha)})} \biggr].
\]
\end{theorem}
%
%
\begin{remark}
Condition~(\ref{eqdecayrate2}) is only required for
Lemma~\ref{thmmaximumdeviation}. As commented by \citet{wu2011},
it can be reduced to $\alpha p>2$ for linear processes. See
Remark~2 of their paper for more details.
\end{remark}

The key step for proving Theorem~\ref{thmthresholding} is to
establish a~convergence rate for the maximum deviation of sample
autocovariances. The following lemma is adapted from Theorem 3 of
\citet{wu2011}, where the asymptotic distribution of the maximum
deviation was also studied.\vadjust{\goodbreak}
%
%
\begin{lemma}
\label{thmmaximumdeviation}
Assume the conditions of Theorem~\ref{thmthresholding}. Then
\[
\lim_{\T\rightarrow\infty} P\Biggl({\max_{1\leq k<\T}}
|\hat\gamma_k-\E\hat\gamma_k|
\leq c_p'\sqrt{\frac{\log T}{T}}\Biggr) =1,
\]
%
where $c_p'=6(p+4) e^{p/4} \|X_0\|_4 \Theta_4$.
\end{lemma}
\begin{pf*}{Proof of Theorem~\ref{thmthresholding}}
Let $B_{\T}=\lfloor(T/\log T)^{{1}/({2(1+\alpha)})} \rfloor$, and
$\hat\Sigma_{\T,B_{\T}}$ be the banded sample covariance matrix
with the
rectangular kernel. Recall that $b_{\T}=({2}/{T})\sum_{k=1}^{B_{\T}}
k|\gamma_k| + 2\sum_{k=B_{\T}+1}^{\T-1}|\gamma_k|$ from
(\ref{eqbias}). By Lemma~\ref{thmmaximumdeviation}, we have
%
\begin{equation}
\label{eq12}
\lambda({\hat\Sigma_{\T,B_{\T}} - \Sigma_{\T}})
= O_P\Biggl( B_{\T}\sqrt{\frac{\log T}{T}} + b_{\T}\Biggr).
\end{equation}
Write the thresholded estimate $\hat\Gamma_{\T,A_{\T}} =
\hat\Gamma_{\T, A_{\T},1} + \hat\Gamma_{\T, A_{\T},2}$, where
\[
\hat\Gamma_{\T,A_{\T},1} = \bigl( \hat\gamma_{s-t}
\hh{1}_{|\hat\gamma_{s-t}|\geq A_{\T}, |s-t|\leq B_{\T}}
\bigr)_{1\leq s,t\leq\T}
\]
and
\[
\hat\Gamma_{\T,A_{\T},2} = \bigl( \hat\gamma_{s-t}
\hh{1}_{|\hat\gamma_{s-t}|\geq A_{\T},
|s-t|> B_{\T}}\bigr)_{1\leq s,t\leq\T}.
\]
By Ger\v{s}gorin's theorem, it is easily seen that
%
\begin{equation}
\label{eq14}
\lambda({\hat\Gamma_{\T,A_{\T},1}-\hat\Sigma_{\T,B_{\T}}})
\leq A_{\T} B_{\T}
= O\Biggl(B_{\T}\sqrt{\frac{\log T}{T}}\Biggr).
\end{equation}
On the other hand,
\begin{eqnarray*}
\lambda({\hat\Gamma_{\T,A_{\T},2}}) &\leq& 2 \Biggl(
\sum_{k=B_{\T}+1}^{\T} |\hat\gamma_k
- \E\hat\gamma_k
|\hh{1}_{|\gamma_k|<A_{\T}/2,|\hat\gamma_k|\geq A_{\T}}
\\
&&\hspace*{10pt}{} + \sum_{k=B_{\T}+1}^{\T} |\hat\gamma_k
- \E\hat\gamma_k|
\hh{1}_{|\gamma_k|\geq A_{\T}/2,|\hat\gamma_k|\geq A_{\T}}
+ \sum_{k=B_{\T}+1}^{\T} |\E\hat\gamma_k|\Biggr) \\
&=&\!: 2(I_{\T}+\II_{\T}+\III_{\T}).
\end{eqnarray*}
The term $\III_{\T}$ is dominated by $b_{\T}$. By
Lemma~\ref{thmmaximumdeviation}, we know
%
\begin{equation}
\label{eq15}
\lim_{\T\rightarrow\infty}P(I_{\T} >0)
\leq\lim_{\T\rightarrow\infty}
P\Bigl(\max_{1\leq k \leq\T-1} |\hat\gamma_k
- \E\hat\gamma_k| \geq A_{\T}/2 \Bigr)
=0.
\end{equation}
For the remaining term $\II_{\T}$, note that the number of $\gamma_k$
such that $k>B_{\T}$ and $|\gamma_k|\geq A_{\T}/2$ is bounded by
$2\sum_{k=B_{\T}+1}^{\T} |\gamma_k|/A_{\T}$; therefore by
Lemma~\ref{thmmaximumdeviation}
%
\begin{equation}
\label{eq16}
\II_{\T} \leq C ({B_{\T}^{-\alpha}}/{A_{\T}})
\cdot\max_{1\leq k\leq\T-1} |\hat\gamma_k-\E\hat\gamma_k|
= O_P(B_{\T}^{-\alpha}).
\end{equation}
Putting~(\ref{eq12}),~(\ref{eq14}),~(\ref{eq15}) and
(\ref{eq16}) together, the proof is complete.
\end{pf*}
%
%
\begin{remark}
If the mean $\E X_1$ is unknown, we need to replace $\hat\gamma_k$
by $\hat\gamma_k^c$ (Remark~\ref{rkcentering}). Since
Lemma~\ref{thmmaximumdeviation} still holds when $\hat\gamma_k$
are replaced by $\hat\gamma_k^c$ [\citet{wu2011}],
Theorem~\ref{thmthresholding} remains true for $\hat\gamma_k^c$.
\end{remark}

\section{A~simulation study}
\label{secsimulation}

The thresholded estimate is desirable in that it can lead to a~better estimate when there are a~lot of zeros or very weak
autocovariances. Unfortunately, due to technical difficulties, the
theoretical result (cf. Theorem~\ref{thmthresholding}) does not
reflect this advantage. We show by simulations that thresholding
does have a~better finite sample performance over banding when the
true autocovariance matrix is sparse.

Consider two linear processes $X_t=\sum_{s=0}^\infty a_s
\varepsilon_{t-s}$ and $Y_t=\sum_{s=0}^\infty b_s\varepsilon_{t-s}$,
where $a_0=b_0=1$, and when $s>0$
\[
a_s = c s^{-(1+\alpha)},\qquad
b_s=c (s/2)^{-(1+\alpha)}\hh{1}_{s
\ \mathrm{is}\  \mathrm{even}}
\]
for some $c>0$ and $\alpha>0$; and $\varepsilon_s$'s are taken as i.i.d.
$\mathcal{N}(0,1)$. Let $\gamma^{X}_k$, $\Sigma^X_{\T}$, and
$\gamma^{Y}_{k}$, $\Sigma^Y_{\T}$ denote the autocovariances and
autocovariance matrices of the two processes, respectively. It is
easily seen that $\gamma^{Y}_k=0$ if $k$ is odd. In fact, $(Y_t)$
can be obtained by interlacing two i.i.d. copies of $(X_t)$. For a~given set of centered observations $X_1, X_2, \ldots, X_{\T}$,
assuming that its true autocovariance matrix is known, for a~fair
comparison we choose the optimal bandwidth $B_{\T}$ and threshold
$A_{\T}$ as
\[
\hat A_{\T}^X = \mathop{\arg\min}_{l\in\{|\hat\gamma^X_0|,
|\hat\gamma^X_1|,\ldots,|\hat\gamma^X_{\T-1}|\}}
\lambda({\hat\Gamma^X_{\T,l} - \Sigma^X_{\T}}),\qquad
\hat B_{\T}^X = \mathop{\arg\min}_{1 \leq k \leq\T}
\lambda({\hat\Sigma^X_{\T,k} - \Sigma^X_{\T}}).
\]
The two parameters for the $(Y_t)$ process are chosen in the same
way. In all the simulations we set $c = 0.5$. For different
combinations of the sample size~$T$ and the parameter $\alpha$ which
controls the decay rate of autocovariances, we report the average
distances in term of the operator norm of the two estimates
$\hat\Sigma_{\T,\hat B_{\T}}$ and $\hat\Gamma_{\T,\hat A_{\T}}$ from
$\Sigma_{\T}$, as well as the standard errors based on 1000
repetitions. We also give the average bandwidth of
$\hat\Sigma_{\T,\hat B_{\T}}$. Instead of reporting the average
threshold for $\hat\Gamma_{\T,\hat A_{\T}}$, we provide the average
number of nonzero autocovariances appearing in the estimates, which is
comparable to the average bandwidth of $\hat\Sigma_{\T,\hat
B_{\T}}$.\vspace*{1pt}

\begin{table}
\caption{Errors under operator norm for $(X_t)$}
\label{tabx}
\begin{tabular*}{\tablewidth}{@{\extracolsep{\fill}}lc
ccd{2.8}d{2.8}d{2.8}@{\hspace*{-2pt}}}
\hline
& \multicolumn{2}{c}{$\bolds{T=100}$} &\multicolumn{2}{c}{$\bolds{T=200}$}
& \multicolumn{2}{c@{}}{$\bolds{T=500}$}
\\[-4pt]
& \multicolumn{2}{c}{\hrulefill} &\multicolumn{2}{l}{\hrulefill\hspace*{2pt}}
& \multicolumn{2}{c@{}}{\hrulefill}
\\
& \multicolumn{1}{c}{\textbf{Error}} & \multicolumn{1}{c}{\textbf{BW}}
& \multicolumn{1}{c}{\textbf{Error}} & \multicolumn{1}{c}{\textbf{BW}}
& \multicolumn{1}{c}{\textbf{Error}} & \multicolumn{1}{c@{}}{\textbf{BW}} \\
\hline
0.2 & 2.94 (1.17) & 9.55 (6.60) & 3.01 (1.22) & 13.4\ (7.67)& 2.96\
(1.23) & 23.4\ (13.1)\\
& 3.66 (1.07) & 5.40 (4.87) & 3.88 (1.14) & 7.39\ (5.81) & 4.08\ (1.17)
& 12.5\ (10.1)\\
& 6.98 (2.63) & & 8.12 (2.85) & & 10.57\ (3.93) & \\
[3pt]
0.5 & 1.52 (0.68) & 6.31 (4.58) & 1.38 (0.60) & 8.46\ (5.57) & 1.15\
(0.50) & 11.9\ (7.67)\\
& 1.90 (0.64) & 3.49 (2.56) & 1.89 (0.59) & 4.15\ (3.07) & 1.74\ (0.54) &
5.15\ (3.27)\\
& 5.55 (2.37) & & 6.73 (2.91) & & 8.88\ (3.28) & \\
[3pt]
$1$ & 0.82 (0.39) & 4.04 (2.33) & 0.69 (0.32) & 4.62\ (2.47) & 0.52\ (0.24) &
5.68\ (3.06)\\
& 1.03 (0.38) & 2.24 (0.87) & 0.95 (0.32) & 2.29\ (0.74) & 0.81\ (0.29) &
2.58\
(0.83)\\
& 4.80 (2.14) & & 6.05 (2.25) & & 7.81\ (2.64) & \\
\hline
\end{tabular*}
\legend{``Error'' refers to the average distance
between the estimates and the true autocovariance matrix under
the operator norm, and ``BW'' refers to the average bandwidth of
the banded estimates, and the average number of nonzero
sub-diagonals (including the diagonal) for the thresholded
ones. The numbers 0.2, 0.5 and 1 in the first column are values of
$\alpha$. For each combination of $T$ and $\alpha$, three lines
are reported, corresponding to banded estimates, thresholded
ones and sample autocovariance matrices, respectively. Numbers in
parentheses are standard errors.}
\end{table}

From Table~\ref{tabx}, we see that for the process $(X_t)$, the
banding method outperforms the thresholding one, but the latter
does give sparser estimates. For the process $(Y_t)$, according to
Table~\ref{taby}, we find that thresholding performs better than
banding when the sample size is not very large ($T=100,200$), and
yields sparser estimates as well. The advantage of thresholding in
error disappears when the sample size is~500. Intuitively
speaking, banding is a~way to threshold according to the truth
(autocovariances with large lags are small), while thresholding is
a~way to threshold according to the data. Therefore, if the
autocovariances are nonincreasing as for the process $(X_t)$, or
if the sample size is large enough, banding is preferable. If the
autocovariances\vadjust{\goodbreak} do not vary regularly as for the process $(Y_t)$
and the sample size is moderate, thresholding is more adaptive. As
a~combination, in practice we can use a~thresholding-after-banding
estimate which enjoys both advantages.

\begin{table}[b]
\caption{Error under operator norm for $(Y_t)$}
\label{taby}
\begin{tabular*}{\tablewidth}{@{\extracolsep{\fill}}lc
ccd{2.8}d{2.8}d{2.8}@{\hspace*{-2pt}}}
\hline
& \multicolumn{2}{c}{$\bolds{T=100}$} &\multicolumn{2}{c}{$\bolds{T=200}$}
& \multicolumn{2}{c@{}}{$\bolds{T=500}$}
\\[-4pt]
& \multicolumn{2}{c}{\hrulefill} &\multicolumn{2}{l}{\hrulefill\hspace*{2pt}}
& \multicolumn{2}{c@{}}{\hrulefill}
\\
& \multicolumn{1}{c}{\textbf{Error}} & \multicolumn{1}{c}{\textbf{BW}}
& \multicolumn{1}{c}{\textbf{Error}} & \multicolumn{1}{c}{\textbf{BW}}
& \multicolumn{1}{c}{\textbf{Error}} & \multicolumn{1}{c@{}}{\textbf{BW}} \\
\hline
0.2 & 3.33 (0.86) & 9.87 (6.89) & 3.54 (0.95) & 13.7\ (7.67) & 3.61\
(1.07) & 24.7\ (13.1)\\
& 3.15 (0.93) & 3.95 (3.50) & 3.43 (1.00) & 5.69\ (4.72) & 3.75\ (1.08)
& 9.23\ (8.04)\\
& 7.21 (4.28) & & 8.69 (4.79) & & 11.1\ (5.31) & \\
[3pt]
0.5 & 1.98 (0.61) & 7.26 (5.32) & 1.88 (0.59) & 9.95\ (6.44) & 1.63\
(0.53) & 16.3\ (10.1)\\
& 1.81 (0.60) & 2.93 (2.41) & 1.81 (0.59) & 3.44\ (2.22) & 1.71\ (0.54) &
4.64\ (2.97)\\
& 5.88 (3.27) & & 7.25 (3.59) & & 9.25\ (3.72) & \\
[3pt]
$1$ & 1.19 (0.41) & 5.31 (3.33) & 1.01 (0.35) & 6.20\ (3.58) & 0.79\ (0.28)
& 8.28\ (4.95)\\
& 1.02 (0.39) & 2.16 (0.65) & 0.92 (0.32) & 2.21\ (0.57) & 0.80\ (0.28)
& 2.52\ (0.77)\\
& 5.09 (2.77) & & 6.39 (2.79) & & 8.18\ (2.91) & \\
\hline
\end{tabular*}
\end{table}

Apparently our simulation is a~very limited one, because we assume
that the true autocovariance matrices are known. Practitioners would
need a~method to choose the bandwidth and/or threshold from the
data. Although theoretical results suggest convergence rates of
banding and thresholding parameters which lead to optimal convergence
rates of the estimates, they do not offer much help for finite
samples. The issue was addressed by \citet{wuss2009}
incorporating the
idea of risk minimization from \citet{bickel2008a} and the technique
of subsampling from \citet{politis1999}, and by \citet{mcmurry2010}
using the rule introduced in \citet{politis2003} for selecting the
bandwidth in spectral density estimation. An alternative method is to
use the block length selection procedure in \citet{buhlmann1999} which
is designed for spectral density estimation. We shall study other
data-driven methods in the future.

\section{Moment inequalities}\label{sec6}
This section presents some moment inequalities that will be useful in
the subsequent proofs. In Lemma~\ref{lemmmtinq}, the case $1 < p \leq
2$ follows from \citet{burkholder1988} and the other case $p > 2$ is
due to \citet{rio2009}. Lemma~\ref{thmfacts} is adopted from
Proposition 1 of \citet{wu2011}.
%
%
\begin{lemma}[{[\citet{burkholder1988}, \citet{rio2009}]}]
\label{lemmmtinq}
Let $p > 1$ and $p' = \min\{p, 2\}$; let~$D_t$, $1\le t \le T$, be
martingale differences, and $D_t \in\LLL^p$ for every $t$. Write
$M_{\T} = \sum_{t=1}^{\T} D_t$. Then
%
\begin{equation}
\label{eqmarzyg}
\|M_{\T}\|_p^{p'} \leq\CC_p^{p'} \sum_{t=1}^{\T} \|D_t\|_p^{p'}
\qquad\mbox{where }
\mathcal{C}_p =\cases{
(p-1)^{-1}, &\quad if $1<p\leq2$,\vspace*{2pt}\cr
\sqrt{p-1}, &\quad if $p > 2$.}\hspace*{-22pt}
\end{equation}
\end{lemma}

It is convenient to use $m$-dependence approximation for processes
with the form~(\ref{eqcausal}). For $t\in\Z$, define $\FF_t =
\langle\varepsilon_t, \varepsilon_{t+1}, \ldots\rangle$ be the
$\sigma$-field generated by the innovations $\varepsilon_t,
\varepsilon_{t+1}, \ldots,$ and the projection operator $\HH_t
(\cdot) = \E(\cdot|\FF_t)$ and $\PP_t(\cdot) = \HH_t(\cdot) -
\HH_{t+1}(\cdot)$. Observe that $(\PP_{-t}(\cdot))_{t \in\Z}$ is
a~martingale difference sequence with respect to the filtration
$(\FF_{-t})$. For $m\geq0$, define $\tilde X_t=\HH_{t-m} X_t$;
then $\|X_t-\tilde X_t\|_p \leq\CC_p \Psi_p(m+1)$ [see
Proposition 1 of \citet{wu2011} for a~proof], and $(\tilde
X_t)_{t\in\Z}$ is an $(m+1)$-dependent sequence.
%
%
\begin{lemma}
\label{thmfacts}
Assume\vspace*{1pt} $\E X_t=0$ and $p > 1$. For $m \geq0$, define $\tilde X_t
= \HH_{t-m} X_t = \E(X_t | \FF_{t-m})$. Let $\tilde
\delta_p(\cdot)$ be the physical dependence measure for the
sequence~$(\tilde X_t)$. Then
%
\begin{eqnarray}
\label{eqfact1}
\|\PP_0X_t\|_p &\leq&\delta_p(t)\quad
\mbox{and}\quad \tilde\delta_{p}(t) \leq\delta_p(t),
\\
\label{eqfact4}
|\gamma_k| &\leq&\zeta_2(k)\qquad \mbox{where }
\zeta_p(k):=\sum_{j=0}^\infty
\delta_p(j)\delta_p(j+k),\\[-18pt]\nonumber
\end{eqnarray}
\begin{eqnarray}
\label{eqfact5}
\Biggl\|\sum_{s,t=1}^{\T} c_{s,t} (X_sX_t-\gamma_{s-t})\Biggr\|_{p/2}
&\leq&4\CC_{p/2}\CC_p \Theta_p^2 \mathcal{B}_{\T} \sqrt{T}
\qquad\mbox{when } p \geq4, \\
\label{eqfact8}
\Biggl\|\sum_{t=1}^{\T} c_t(X_t-\tilde X_t)\Biggr\|_p
&\leq&\mathcal{C}_p \mathcal{A}_{\T} \Theta_{p}(m+1)
\qquad\mbox{when } p \geq2,
\end{eqnarray}
where
\[
\mathcal{A}_{\T} = \Biggl(\sum_{t=1}^{\T} |c_t|^2\Biggr)^{1/2}
\quad\mbox{and}\quad
\mathcal{B}_{\T}^2 =\max\Biggl\{\max_{1\leq t \leq\T} \sum
_{s=1}^{\T} c_{s,t}^2,
\max_{1\leq s \leq\T} \sum_{t=1}^{\T} c_{s,t}^2\Biggr\}.
\]
\end{lemma}

\section{Large deviations for quadratic forms}
\label{secld}

In this section we prove a~result on probabilities of large deviations
of quadratic forms of stationary processes, which take the form
\[
Q_{\T} = \sum_{1 \leq s \leq t \leq\T} a_{s,t}X_sX_t.
\]
The coefficients $a_{s,t}=a_{\T,s,t}$ may depend on $T$, but we
suppress $T$ from subscripts for notational simplicity. Throughout
this section we assume that ${\sup_{s,t}} |a_{s,t}| \leq1$, and
$a_{s,t}=0$ when $|s-t|>B_{\T}$, where $B_{\T}$ satisfies $B_{\T}
\rightarrow\infty$, and $B_{\T} = O(T^{\gamma})$ for some
$0<\gamma<1$.

Large deviations for quadratic forms of stationary processes have
been extensively studied in the literature. \citet{bryc1997} and
\citet{bercu1997} obtained the \textit{large deviation principle}
[\citet{dembo1998}] for Gaussian processes. \citet{gamboa1999}
considered the functional large deviation principle.
\citet{bercu2000} obtained a~more accurate expansion of the tail
probabilities. \citet{zani2002} extended the results of
\citet{bercu1997} to locally stationary Gaussian processes. In
fact, our result is more relevant to the so-called \textit{moderate
deviations} according to the terminology of \citet{dembo1998}.
\citet{bryc1997} and \citet{kakizawa2007} obtained
\textit{moderate deviation principles} for quadratic forms of Gaussian processes.
\citet{djellout2006} studied moderate deviations of periodograms
of linear processes. \citet{bentkus1976} considered the
Cram\'er-type moderate deviation for spectral density estimates of
Gaussian processes; see also \citet{saulis1991}. \citet{liu2010}
derived the Cram\'er-type moderate deviation for maxima of
periodograms under the assumption that the process consists of i.i.d.
random variables.

For our purpose, on one hand, we do not need a~result that is as
precise as the moderate deviation principle or the Cram\'er-type
moderate deviation. On the other hand, we need an upper bound for
the tail probability under less restrictive conditions.
Specifically, we would like to relax the Gaussian, linear or i.i.d.
assumptions which were made in the precedent works.
\citet{rudzkis1978} provided a~result in this fashion under the
assumption of boundedness of the cumulant spectra up to a~finite
order. While this type of assumption holds under certain mixing
conditions, the latter themselves are not easy to verify in
general and many well-known examples are not strong mixing
[\citet{andrews1984}]. We mean to impose alternative conditions
through physical dependence measures, which are easy to use in
many applications [\citet{wu2005}]. Furthermore, our result can be
sharper; see Remark~\ref{rkrudzkis}.\vadjust{\goodbreak}

Our main tool is the $m$-dependence approximation. In the next
lemma we use dependence measures to bound the $\LLL^p$ norm of the
distance between~$Q_{\T}$ and the $m$-dependent version $\tilde Q_{\T}$.
The proof and a~few remarks on the optimality of the result are
given in the supplementary article [\citet{wu2011s}].
%
%
\begin{lemma}
\label{thmmappquadratic}
Assume $X_t \in\LLL^p$ with $p\geq4$, $\E X_t=0$ and
$\Theta_{p}<\infty$. Let $\tilde X_t = \HH_{t-m_{\T}}X_t$ and
$\tilde Q_{\T}
= \sum_{1\leq s\leq t \leq\T} a_{s,t} \tilde X_s \tilde X_t$; then
\begin{eqnarray*}
&&
\|\E_0 Q_{\T} - \E_0 \tilde Q_{\T}\|_{p/2}\\
&&\qquad\leq 4 \Theta_p(m_{\T})^2 + 11 (p-2) \Theta_p \sqrt{TB_{\T}}
\Theta_p(m_{\T}) \\
&&\qquad\quad{} +(p-2) \sqrt{TB_{\T}} \bigl[ 3\Theta_p(\lfloor m_{\T}/2 \rfloor
)\Delta_p(m_{\T})
+ 3\Theta_p(m_{\T})\Delta_p(\lfloor m_{\T}/2 \rfloor)\bigr].
\end{eqnarray*}
\end{lemma}

The following theorem is the main result of this section.
%
%
\begin{theorem}
\label{thmldquadratic}
Assume $X_t \in\LLL^p$, $p>4$, $\E X_t=0$, and
$\Theta_p(m)=O(m^{-\alpha})$. Set $c_p = (p+4)e^{p/4}\Theta_4^2$.
For any $M>1$, let $x_{\T} = 2c_p \sqrt{TMB_{\T}\log B_{\T}}$. Assume
that $B_{\T}\rightarrow\infty$ and $B_{\T}=O(T^\gamma)$ for some
$0<\gamma<1$. Then for any $\gamma< \beta< 1$, there exists a~constant $C_{p,M,\beta}>0$ such that
\begin{eqnarray*}
&&
P (|\E_0 Q_{\T}| \geq x_{\T} )\\
&&\qquad\leq C_{p,M,\beta} x_{\T}^{-p/2} (\log T)
[(TB_{\T})^{p/4}T^{-\alpha\beta p/2} + T B_{\T}^{p/2-1-\alpha\beta p/2} + T]\\
&&\qquad\quad{} + C_{p,M,\beta} B_{\T}^{-M}.
\end{eqnarray*}
\end{theorem}
%
%
\begin{remark}
\label{rkrudzkis}
\citet{rudzkis1978} proved that if $p=4k$ for some $k\in\N$, then
\[
P (|\E_0 Q_{\T}| \geq x_{\T} )
\leq C x_{\T}^{-p/2}(TB_{\T})^{p/4},
\]
which can be obtained by using Markov inequality and
(\ref{eqfact5}) under our framework. The upper bound given in
Theorem~\ref{thmldquadratic} has a~smaller order of magnitude. We
note that \citet{rudzkis1978} also proved a~stronger exponential
inequality under strong moment conditions. They required the
existence of every moment and the absolute summability of cumulants
of every order.
\end{remark}
\begin{pf*}{Proof of Theorem~\ref{thmldquadratic}}
Without loss of generality, assume $B_{\T} \le T^{\gamma}$. For
$\gamma<\beta<1$, let
$m_{\T}=\lfloor T^{\sqrt{\beta}} \rfloor$, $\tilde X_t = \HH
_{t-m_{\T
}}X_t$ and
\[
\tilde Q_{\T} = \sum_{1\leq s\leq t\leq\T}
a_{s,t} \tilde X_s \tilde X_{t}.
\]
By Lemma~\ref{thmmappquadratic} and~(\ref{eqfact5}), we have
%
\begin{eqnarray}
\label{eq2}
&&P\bigl[|\E_0(Q_{\T}-\tilde Q_{\T})| \geq c_p M^{1/2}\sqrt{TB_{\T
}(\log B_{\T})} \bigr]\nonumber\\[-8pt]\\[-8pt]
&&\qquad \leq C_{p,M} x_{\T}^{-p/2} (TB_{\T})^{p/4} T^{-\alpha
\sqrt{\beta} p/2}.
\nonumber
\end{eqnarray}
Split $[1,T]$ into blocks $\mathcal{B}_1,\ldots,\mathcal{B}_{b_{\T}}$
of size $2m_{\T}$, and define
\[
Q_{\T,k} = \sum_{t \in\mathcal{B}_k} \sum_{1\leq s \leq t}
a_{s,t}\tilde X_s \tilde X_t.
\]
By Corollary 1.7 of \citet{nagaev1979} and~(\ref{eqfact5}), we know
for any $M>1$, there exists a~constant $C_{p,M,\beta}$ such that
%
\begin{eqnarray}
\label{eq3}
&&P \bigl[|\E_0 \tilde Q_{\T}|
\geq c_p \sqrt{TMB_{\T}(\log B_{\T})}\bigr]\nonumber\\[-2pt]
&&\qquad\leq\sum_{k=1}^{b_{\T}} P \biggl(|\E_0 Q_{\T,k}|
\geq\frac{x_{\T}}{C_{p,M,\beta}} \biggr)  + \biggl[\frac{C_{p,M,\beta} Tm_{\T}^{-1} (m_{\T}B_{\T})^{p/4}}
{(TB_{\T})^{p/4}}\biggr]^{C_{p,M,\beta}}
\nonumber\\[-9pt]\\[-9pt]
&&\qquad\quad{}+ C_{\beta} \exp\biggl\{\frac{c_p^2 (\log B_{\T})}
{(p+4)^2 e^{p/2} \Theta_4^4}\biggr\} \nonumber\\[-2pt]
&&\qquad\leq \sum_{k=1}^{b_{\T}} P (|\E_0 Q_{\T,k}|
\geq x_{\T}/C_{p,M,\beta} )
+ C_{p,M,\beta} (B_{\T}^{-M} + T^{-M}).\nonumber
\end{eqnarray}
By Lemma~\ref{thmldtntermediate}, we have
%
\begin{eqnarray}
\label{eq4}
&&P (|\E_0 Q_{\T,k}| \geq x_{\T}/C_{p,M,\beta} )\nonumber\\[-2pt]
&&\qquad\leq C_{p,M,\beta} x_{\T}^{-p/2} (\log T) \\[-2pt]
&&\qquad\quad{} \times\bigl[ \bigl(T^{\sqrt{\beta}} B_{\T}\bigr)^{p/4}T^{-\alpha
\beta p/2}
+ T^{\sqrt{\beta}}B_{\T}^{p/2-1-\alpha\beta p/2}
+ T^{\sqrt{\beta}}\bigr].
\nonumber
\end{eqnarray}
Combining~(\ref{eq2}),~(\ref{eq3}) and~(\ref{eq4}), the proof is
complete.
\end{pf*}
%
%
\begin{lemma}
\label{thmldtntermediate}
Assume $X_t \in\LLL^p$ with $p>4$, $\E X_t=0$, and
$\Theta_p(m)=O(m^{-\alpha})$. If $x_{\T}>0$ satisfies
$T^{\delta}\sqrt{TB_{\T}}=o(x_{\T})$ for some $\delta>0$, then for any
$0<\beta<1$, there exists a~constant $C_{p,\delta,\beta}$ such
that
\begin{eqnarray*}
P(|\E_0 Q_{\T}| \geq x_{\T})
&\leq& C_{p,\delta,\beta} x_{\T}^{-p/2} (\log T) \\[-2pt]
&&{} \times[(TB_{\T})^{p/4}T^{-\alpha\beta p/2}
+ TB_{\T}^{p/2-1-\alpha\beta p/2} + T].
\end{eqnarray*}
\end{lemma}
\begin{pf}
For $j\geq1$, define $m_{\T,j}=\lfloor T^{\beta^j} \rfloor$,
$X_{t,j}=\HH_{t-m_{\T,j}}X_t$ and
\[
Q_{\T,j} = \sum_{1\leq s\leq t\leq\T} a_{s,t}X_{s,j}X_{t,j}.
\]
Let $j_{\T}=\lceil-\log(\log T)/(\log\beta) \rceil$. Note that
$m_{\T,j_{\T}} \leq e$. By Lemma~\ref{thmmappquadratic} and
(\ref{eqfact5}),
%
\begin{equation}
\label{eq1}
P[|\E_0(Q_{\T}-Q_{\T,1})| \geq x_{\T}/j_{\T} ]
\leq C_{p,\beta} (\log T)^{1/2} x_{\T}^{-p/2} (TB_{\T})^{p/4}
T^{-\alpha\beta p/2}.\hspace*{-22pt}
\end{equation}
Let $j_{\T}'$ be the smallest $j$ such that $m_{\T,j}< B_{\T}/4$. For
$1\leq j < j_{\T}'$, split $[1,T]$ into blocks $\mathcal{B}^{(j)}_1,
\ldots, \mathcal{B}^{(j)}_{b_{\T,j}}$ of size $B_{\T}+m_{\T,j}$. Define
\[
R_{\T,j,b}=\sum_{t\in\mathcal{B}^{(j)}_b}\sum_{1 \leq s \leq t}
a_{s,t} X_{s,j}X_{t,j} \quad\mbox{and}\quad
R_{\T,j,b}'=\sum_{t\in\mathcal{B}^{(j)}_b}\sum_{1 \leq s \leq t}
a_{s,t} X_{s,j+1}X_{t,j+1}.\vadjust{\goodbreak}
\]
By Corollary 1.6 of \citet{nagaev1979} and~(\ref{eqfact5}), we have
for any $C>2$
%
\begin{eqnarray}
\label{eq17}\qquad
P \biggl[|\E_0(Q_{\T,j} - Q_{\T,j+1})|>\frac{x_{\T}}{2j_{\T
}}\biggr]
&\leq& \sum_{b=1}^{b_{\T,j}}
P \biggl[|\E_0(R_{\T,j,b}-R_{\T,j,b}')|
\geq\frac{x_{\T}}{C j_{\T}}\biggr] \\
\label{eq18}
&&{} + 2\biggl[\frac{64 C e^2\Theta_4^4TB_{\T}j_{\T}^2}{x_{\T
}^2}\biggr]^{C/4}.
\end{eqnarray}
It is clear that for any $M>1$, there exists a~constant
$C_{M,\delta,\beta}$ such that the term in~(\ref{eq18}) is less
than $C_{M,\delta,\beta} x_{\T}^{-M}$. For~(\ref{eq17}), by
Lemma~\ref{thmmappquadratic} and~(\ref{eqfact5})
\begin{eqnarray*}
&&\sum_{b=1}^{b_{\T,j}}
P\biggl[|\E_0(R_{\T,j,b}-R_{\T,j,b}')|
\geq\frac{x_{\T}}{C j_{\T}}\biggr] \\
&&\qquad \leq C_{p,\beta} {T}{(m_{\T,j})^{-1}}
\cdot(\log T)^{1/2}\cdot x_{\T}^{-p/2}
\cdot(m_{\T,j}B_{\T})^{p/4} \cdot m_{\T,j+1}^{-\alpha p/2} \\
&&\qquad \leq C_{p,\beta} x_{\T}^{-p/2}\cdot(\log T)^{1/2}TB_{\T}^{p/4}
\cdot(m_{\T,j})^{p/4-1-\alpha\beta p/2}.
\end{eqnarray*}
Depending on whether the exponent $p/4-1-\alpha\beta p/2$ is
positive or not, the term $(m_{\T,j})^{p/4-1-\alpha\beta p/2}$ is
maximized when $j=1$ or $j=j_{\T}'-1$, respectively, and we have
%
\begin{eqnarray}
\label{eq25}\quad
&&\sum_{b=1}^{b_{\T,j}}
P\biggl[|\E_0(R_{\T,j,b}-R_{\T,j,b}')|
\geq\frac{x_{\T}}{C j_{\T}} \biggr] \nonumber\\[-8pt]\\[-8pt]
&&\qquad \leq C_{p,\beta} x_{\T}^{-p/2}\cdot(\log T)^{1/2} \cdot
[(TB_{\T})^{p/4}T^{-\alpha\beta p/2}
+ T B_{\T}^{p/2-1-\alpha\beta p/2}].
\nonumber
\end{eqnarray}
Combining~(\ref{eq1}),~(\ref{eq17}),~(\ref{eq18}) and
(\ref{eq25}), we have shown that
%
\begin{eqnarray}
\label{eq29}\qquad
&&P (|\E_0 Q_{\T}| \geq x_{\T})\nonumber\\
&&\qquad\leq P(|\E_0 Q_{\T,j_{\T}'}| \geq x_{\T}/2)
+ C_{p,M,\delta,\beta} x_{\T}^{-M} \\
&&\qquad\quad{} + C_{p,M,\delta,\beta} x_{\T}^{-p/2} (\log T)
[(TB_{\T})^{p/4}T^{-\alpha\beta p/2}
+ TB_{\T}^{p/2-1-\alpha\beta p/2}].
\nonumber
\end{eqnarray}

To deal with\vspace*{2pt} the probability concerning $Q_{\T,j_{\T}'}$ in
(\ref{eq29}), we split $[1,T]$ into blocks
$\mathcal{B}_1,\ldots,\mathcal{B}_{b_{\T}}$ with size $2B_{\T}$, and
define the block sums
\[
R_{\T,j_{\T}',b}=\sum_{t\in\mathcal{B}_b}\sum_{1 \leq s \leq t}
a_{s,t} X_{s,j_{\T}'}X_{t,j_{\T}'}. 
\]
Similarly as~(\ref{eq17}) and~(\ref{eq18}), there exists a~constant $C_{p,M,\delta, \beta}>2$ such that
\[
P(|\E_0 Q_{\T,j_{\T}'}| \geq{x_{\T}}/{2})
\leq\sum_{b=1}^{b_{\T}} P\biggl(|\E_0 R_{\T,j_{\T
}',b}|
\geq\frac{x_{\T}}{C_{p,M,\delta,\beta}}\biggr) + C_{p,M,\delta
,\beta} x_{\T}^{-M}.
\]
By Lemma~\ref{thmldnn}, we have
\[
P(|\E_0 R_{\T,j_{\T}',b}|
\geq C_{p,M,\delta,\beta}^{-1} x_{\T}) \leq C_{p,M,\delta
,\beta} x_{\T}^{-p/2} (\log T)
(B_{\T}^{p/2-\alpha\beta p/2} + B_{\T});
\]
and it follows that for some constant $C_{p,\delta,\beta}>0$,
%
\begin{equation}
\label{eq30}
P(|\E_0 Q_{\T,j_{\T}'}| \geq x_{\T}/2) \leq
C_{p,\delta,\beta} x_{\T}^{-p/2} (\log T) T
(B_{\T}^{p/2-1-\alpha\beta p/2}+1).
\end{equation}
The proof is completed by combining~(\ref{eq29}) and
(\ref{eq30}).
\end{pf}

In the next lemma we consider $Q_{\T}$ when the restriction
$a_{s,t}=0$ for $|s-t|>B_{\T}$ is removed. To avoid confusion, we use
a~new symbol. Let
\[
R(T,m)=\sum_{1 \le s \le t \le\T}
c_{s,t} (\HH_{s-m}X_s)(\HH_{t-m}X_t).
\]
For $x_{\T}>0$, define
\[
U(T,m,x_{\T}) = \sup_{\{c_{s,t}\}} P[|\E_0 R({T,m})| \geq x_{\T}],
\]
where the supremum is taken over all arrays $\{c_{s,t}\}$ such
that $|c_{s,t}| \leq1$. We use $R_{\T}$ and $U(T,x_{\T})$ as shorthands
for $R({T,\infty})$ and $U(T,\infty,x_{\T})$, respectively.
%
%
\begin{lemma}
\label{thmldnn}
Assume $X_t \in\LLL^p$ with $p>4$, $\E X_t=0$, and
$\Theta_p(m)=O(m^{-\alpha})$. If $x_{\T}>0$ satisfies
$T^{1+\delta}=o(x_{\T})$ for some $\delta>0$, then for any $0<\beta<1$,
there exists a~constant $C_{p,\delta,\beta}$ such that
\[
P (|\E_0 R_{\T}|\geq x_{\T} )
\leq C_{p,\delta,\beta}
x_{\T}^{-p/2} (\log T)(T^{p/2-\alpha\beta p/2}+T).
\]
%
\end{lemma}
\begin{pf}
Let $m_{\T}=\lfloor T^{\beta} \rfloor$ and $\tilde R_{\T
}:=R(T,m_{\T})$.
By Lemma~\ref{thmmappquadratic} and~(\ref{eqfact5}),
\[
P[|\E_0(R_{\T}-\tilde R_{\T})| \geq x_{\T}/2 ]
\leq C_{p} x_{\T}^{-p/2} T^{p/2-\alpha\beta p/2}.
\]
We claim that there exists a~constant
$C_{p,\delta,\beta}$ such that
%
\begin{equation}
\label{eq28}
U(T,m_{\T},x_{\T}/2) \leq C_{p,\delta,\beta}
x_{\T}^{-p/2} (T \log T)(m_{\T}^{p/2-1-\alpha\beta p/2}+1).
\end{equation}
Therefore, the proof is complete by using
\[
P (|\E_0 R_{\T}|\geq x_{\T} )
\leq P[|\E_0(R_{\T}-\tilde R_{\T})| \geq x_{\T}/2 ]
+ U(T,m_{\T},x_{\T}/2).
\]

We need to prove the claim~(\ref{eq28}). Let $z_{\T}$ satisfy
$T^{1+\delta}=o(z_{\T})$. Let $j_{\T}=\lceil-{\log}(\log T)/(\log
\beta) \rceil$, and note that $T^{\beta^{j_{\T}}} \leq e$. Set
$y_{\T} =
z_{\T} / (2j_{\T})$. We consider $U(T,m,z_{\T})$ for an arbitrary
$1<m<T/4$. Set $X_{t,1} := \HH_{t-m}X_t$ and $X_{t,2} :=
\HH_{t-\lfloor m^\beta\rfloor} X_t$. Define
\[
Y_{t,1}=\sum_{s=1}^{t-3m-1}c_{s,t}X_{s,1}
\quad\mbox{and}\quad
Z_{t,1}=\sum_{s=1\vee(t-3m)}^{t}c_{s,t}X_{s,1}
\]
and $Y_{t,2}$, $Z_{t,2}$ similarly by replacing $X_{s,1}$ with
$X_{s,2}$. Observe\vspace*{1pt} that $X_{t,k}$ and $Y_{t,l}$ are independent
for $k,l=1,2$. We first consider $\sum_{t=1}^{\T} (X_{t,1} Z_{t,1} -
X_{t,2}Z_{t,2})$. Split $[1,T]$ into blocks $\mathcal{B}_1, \ldots,
\mathcal{B}_{b_{\T}}$ with size $4 m$, and define $W_{\T,b} =\break
\sum_{t \in\mathcal{B}_b} (X_{t,1}Z_{t,1} - X_{t,2} Z_{t,2})$. Let
$y_{\T}$ satisfy $y_{\T}<z_{\T}/2$ and $T^{1+\delta/2} =
o(y_{\T})$. Since $W_{\T,b}$ and $W_{\T,b'}$ are independent if
$|b-b'| > 1$, by Corollary 1.6 of \citet{nagaev1979},
(\ref{eqfact5}) and Lemma~\ref{thmmappquadratic}, similarly as
(\ref{eq17}) and~(\ref{eq18}), we know for any $M>1$, there exists
a~constant $C_{p,M,\delta,\beta}$ such that
%
\begin{eqnarray}
\label{eq24}
&&P\Biggl[\Biggl|\E_0\Biggl(\sum_{t=1}^{\T}
X_{t,1}Z_{t,1} - X_{t,2}Z_{t,2}\Biggr)\Biggr|
\geq y_{\T} \Biggr] \nonumber\\
&&\qquad \leq C_{p,M,\delta,\beta} y_{\T}^{-M}
+ \sum_{b=1}^{b_{\T}} P(|\E_0 W_{\T,b}|
\geq y_{\T}/C_{M,\delta}) \\
&&\qquad \leq C_{p,M,\delta,\beta} y_{\T}^{-M}
+ C_{p,M,\delta,\beta} y_{\T}^{-p/2}
T m^{p/2-1-\alpha\beta p/2}.\nonumber
\end{eqnarray}
Now we deal with the term $\sum_{t=1}^{\T} (X_{t,1} Y_{t,1} -
X_{t,2} Y_{t,2})$. Split $[1,T]$ into blocks $\mathcal{B}^\ast_1,
\ldots, \mathcal{B}^\ast_{b_{\T}^\ast}$ with size $m$. Define
$R_{\T,b} = \sum_{t \in\mathcal{B}^\ast_b} (X_{t,1}Y_{t,1} -
X_{t,2}Y_{t,2})$. Let $\xi_{b}$ be the $\sigma$-fields generated by
$\{\varepsilon_{l_b}, \varepsilon_{l_b-1}, \ldots\}$, where $l_b =
\max\{\mathcal{B}^\ast_b\}$. Observe that $(R_{\T,b})_{b
\ \mathrm{is}\ \mathrm{odd}}$ is a~martingale sequence with
respect to $(\xi_{b})_{b\ \mathrm{is}\ \mathrm{odd}}$, and so are
$(R_{\T,b})_{b\ \mathrm{is}\ \mathrm{even}}$ and $(\xi_{b})_{b
\ \mathrm{is}\ \mathrm{even}}$. By Lemma 1 of
\citet{haeusler1984} we know for any $M>1$, there exists a~constant
$C_{M,\delta}$ such that
%
\begin{eqnarray}
\label{eq19}
&&P\Biggl[\Biggl|\sum_{t=1}^{\T} (X_{t,1}Y_{t,1}
- X_{t,2}Y_{t,2}) \Biggr| \geq y_{\T} \Biggr] \nonumber\\
&&\qquad \leq C_{M,\delta} y_{\T}^{-M}
+ 4P\Biggl[\sum_{b=1}^{b_{\T}^\ast}
\E(R_{\T,b}^2|\xi_{b-2}) > \frac{y_{\T}^2}
{(\log y_{\T})^{3/2}}\Biggr]\nonumber\\[-8pt]\\[-8pt]
&&\qquad\quad{}+ \sum_{b=1}^{b_{\T}^\ast}
P\Biggl[|R_{\T,b}| \geq\frac{y_{\T}}{\log y_{\T}}\Biggr]\nonumber\\
&&\qquad=: I_{\T} + \II_{\T} + \III_{\T}.\nonumber
\end{eqnarray}
Since $(X_{t,1},X_{t,2})$ and $(Y_{t,1},Y_{t,2})$ are independent,
$R_{\T,b}$ has finite $p$th moment. Using similar arguments as
Lemma~\ref{thmmappquadratic}, we have
\[
\|R_{\T,b}\|_p \leq C_p (mT)^{p/2} m^{-\alpha\beta p};
\]
and it follows that
%
\begin{equation}
\label{eq20}
\III_{\T} \leq C_{p} y_{\T}^{-p}(\log y_{\T})^p T^{p/2+1}
m^{p/2-1-\alpha\beta p}.
\end{equation}
For the second term, let $r_{s-t,k}=\E(X_{s,k}X_{t,k})$ for
$k=1,2$; we have
%
\begin{eqnarray}\label{eq5}\qquad
\sum_{b=1}^{b_{\T}^\ast} \E(R_{\T,b}^2|\xi_{b-2})
& \leq& 2 \sum_{b=1}^{b_{\T}^\ast}
\biggl[\sum_{s,t \in\mathcal{B}_b^\ast}
(r_{s-t,1}Y_{s,1}Y_{t,1}
+ r_{s-t,2}Y_{s,2}Y_{t,2}) \biggr] \nonumber\\[-8pt]\\[-8pt]
& = &\sum_{1 \leq s\leq t \leq\T} a_{s,t,1} X_{s,1}X_{t,1}
+ \sum_{1 \leq s\leq t \leq\T} a_{s,t,2} X_{s,2}X_{t,2}.\nonumber
\end{eqnarray}
By~(\ref{eqfact1}) and~(\ref{eqfact4}), we know $\sum_{l \in\Z}
|r_{l,k}| < \infty$ for $k=1,2$, and hence $|a_{s,t,k}| \leq
CT$. It follows that the expectations of the two terms in
(\ref{eq5}) are all less than $CT^2$, and
%
\begin{equation}
\label{eq21}
\II_{\T} \leq C_{\beta} U\biggl[T,m,
\frac{y_{\T}^2}{T (\log y_{\T})^2}\biggr]
+ C_{\beta} U\biggl[T, \lfloor m^\beta\rfloor,
\frac{y_{\T}^2}{T (\log y_{\T})^2}\biggr].
\end{equation}
Combining~(\ref{eq24}),~(\ref{eq19}),~(\ref{eq20}) and
(\ref{eq21}), we have shown that $U(T,m,z_{\T})$ is bounded from
above by
%
\begin{eqnarray}
\label{eq22}
&& U(T,\lfloor m^\beta\rfloor,z_{\T}-2y_{\T})\nonumber\\
&&\qquad{}+ C_{\beta} U \biggl[T, \lfloor m^\beta\rfloor,
\frac{y_{\T}^2}{T (\log y_{\T})^2}\biggr]
+ C_{\beta} U \biggl[T, m,
\frac{y_{\T}^2}{T (\log y_{\T})^2}\biggr]\nonumber\\[-8pt]\\[-8pt]
&&\qquad{}+ C_{p,M,\delta,\beta} [y_{\T}^{-M}
+ y_{\T}^{-p/2} Tm^{p/2-1-\alpha\beta p/2}\nonumber\\
&&\qquad\hspace*{53pt}{}
+ y_{\T}^{-p}(\log y_{\T})^p T^{p/2+1}m^{p/2-1-\alpha\beta p}].
\nonumber
\end{eqnarray}
Since $\sup_{\{c_{s,t}\}} \|\E_0 R_{\T}\|_{p/2} \leq C_p T$ by
(\ref{eqfact5}), by applying~(\ref{eq22}) recursively to deal with
the last term on the first line of~(\ref{eq22}) for $q$ times such
that $(y_{\T}/T)^{-2^q p}=O [y_{\T}^{-(M+1)} ]$, we have
%
\begin{eqnarray}
\label{eq23}
&&U(T,m,z_{\T})
\leq C_{p,M,\delta,\beta} \bigl[
U(T,\lfloor m^\beta\rfloor,z_{\T}-2y_{\T})
+ y_{\T}^{-p/2} Tm^{p/2-1-\alpha\beta p/2}\hspace*{-25pt}\nonumber\\[-8pt]\\[-8pt]
&&\qquad\hspace*{107pt}{} +
y_{\T}^{-p}(\log y_{\T})^p T^{p/2+1}m^{p/2-1-\alpha\beta p}
+ y_{\T}^{-M}\bigr].\hspace*{-25pt}
\nonumber
\end{eqnarray}
Using the preceding arguments similarly, we can show that when
$1\leq m \leq3$
\[
U[T,m,z_{\T}/(2j_{\T})] \leq C_{M,p,\delta} [ z_{\T}^{-p/2}
(\log T)T
+ z_{\T}^{-p} (\log z_{\T})^{p+1} T^{p/2+1} + z_{\T}^{-M}].
\]
The details of the derivation are omitted. Applying~(\ref{eq23})
recursively for at most $j_{\T}-1$ times, we have the first bound for
$U(T,m,z_{\T})$,
%
\begin{eqnarray}
\label{eq26}
&& U(T,m,z_{\T}) \nonumber\\
&&\qquad \leq C_{p,M,\delta,\beta}^{j_{\T}}
\{U[T,3,z_{\T}/(2j_{\T})]
+ z_{\T}^{-p/2} (\log z_{\T})
T(m^{p/2-1-\alpha\beta p/2}+1)\nonumber\\[-8pt]\\[-8pt]
&&\qquad\quad\hspace*{65.7pt}{}
+ z_{\T}^{-p}(\log z_{\T})^{p+1}
T^{p/2+1}(m^{p/2-1-\alpha\beta p}+1)
+ z_{\T}^{-M}\}\nonumber\\
&&\qquad \leq C_{p,\delta,\beta}^{j_{\T}}(\log z_{\T})^{p+1}
(z_{\T}^{-p/2} T + z_{\T}^{-p} T^{p/2+1})
(m^{p/2-1-\alpha\beta p/2}+1).
\nonumber
\end{eqnarray}
Now plugging~(\ref{eq26}) back into~(\ref{eq22}) for the last two
terms on the first line and using the condition
$T^{1+\delta/2}=o(y_{\T})$, we have
%
\begin{eqnarray}
\label{eq27}
U(T,m,z_{\T} ) &\leq& U(T,\lfloor m^\beta\rfloor,z_{\T} - 2y_{\T}
)\nonumber\\[-8pt]\\[-8pt]
&&{}+ C_{p,\delta,\beta}
[y_{\T}^{-p/2} T (m^{p/2-1-\alpha\beta p/2} + 1)] .\nonumber
\end{eqnarray}
Again by applying~(\ref{eq27}) for at most $j_{\T}-1$ times, we obtain
the second bound for $U(T,m,z_{\T})$:
\[
U(T,m,z_{\T}) \leq C_{p,\delta,\beta}
z_{\T}^{-p/2} (T \log T)(m^{p/2-1-\alpha\beta p/2}+1).
\]
The proof of the claim~(\ref{eq28}) is complete.
\end{pf}

\section{Conclusion}
\label{secconclude}

In this paper we use Toeplitz's connection of eigenvalues of
matrices and Fourier transforms of their entries, and obtain
optimal bounds for tapered covariance matrix estimates by applying
asymptotic results of spectral density estimates. Many problems
are still unsolved; for example, can we improve the convergence
rate of the thresholded estimate in
Theorem~\ref{thmthresholding}? What is the asymptotic
distribution of the maximum eigenvalues of the estimated
covariance matrices? We hope that the approach and results
developed in this paper can be useful for other high-dimensional
covariance matrix estimation problems in time series. Such
problems are relatively less studied compared to the well-known
theory of random matrices which requires i.i.d. entries or multiple
i.i.d. copies.

\section*{Acknowledgments}

We are grateful to an Associate Editor and the referees for their
many helpful comments.

\begin{supplement}[id=suppA]
\stitle{Additional technical proofs}
\slink[doi]{10.1214/11-AOS967SUPP} 
\sdatatype{.pdf}
\sfilename{aos967\_supp.pdf}
\sdescription{We give the proofs of Remark~\ref{rem5} and
Lemma~\ref{thmmappquadratic}, as well as
a~few remarks on Lemma~\ref{thmmappquadratic}.}
\end{supplement}

%

\printaddresses


\begin{thebibliography}{59}

\bibitem[\protect\citeauthoryear{Adenstedt}{1974}]{adenstedt1974}
%
\begin{barticle}[mr]
\bauthor{\bsnm{Adenstedt},~\bfnm{Rolf~K.}\binits{R.~K.}}
(\byear{1974}).
\btitle{On large-sample estimation for the mean of a~stationary random
sequence}.
\bjournal{Ann. Statist.}
\bvolume{2}
\bpages{1095--1107}.
\bid{issn={0090-5364}, mr={0368354}}
\bptok{imsref}%
\end{barticle}
%
\endbibitem

\bibitem[\protect\citeauthoryear{An, Chen and Hannan}{1983}]{an1983}
%
\begin{barticle}[mr]
\bauthor{\bsnm{An},~\bfnm{Hong~Zhi}\binits{H.~Z.}},
\bauthor{\bsnm{Chen},~\bfnm{Zhao~Guo}\binits{Z.~G.}} \AND
\bauthor{\bsnm{Hannan},~\bfnm{E.~J.}\binits{E.~J.}}
(\byear{1983}).
\btitle{The maximum of the periodogram}.
\bjournal{J. Multivariate Anal.}
\bvolume{13}
\bpages{383--400}.
\bid{doi={10.1016/0047-259X(83)90017-9}, issn={0047-259X}, mr={0716931}}
\bptok{imsref}%
\end{barticle}
%
\endbibitem

\bibitem[\protect\citeauthoryear{Andrews}{1984}]{andrews1984}
%
\begin{barticle}[mr]
\bauthor{\bsnm{Andrews},~\bfnm{Donald W.~K.}\binits{D.~W.~K.}}
(\byear{1984}).
\btitle{Nonstrong mixing autoregressive processes}.
\bjournal{J. Appl. Probab.}
\bvolume{21}
\bpages{930--934}.
\bid{issn={0021-9002}, mr={0766830}}
\bptok{imsref}%
\end{barticle}
%
\endbibitem

\bibitem[\protect\citeauthoryear{Bai and Silverstein}{2010}]{bai2010}
%
\begin{bbook}[mr]
\bauthor{\bsnm{Bai},~\bfnm{Zhidong}\binits{Z.}} \AND
\bauthor{\bsnm{Silverstein},~\bfnm{Jack~W.}\binits{J.~W.}}
(\byear{2010}).
\btitle{Spectral Analysis of Large Dimensional Random Matrices},
\bedition{2nd} ed.
\bpublisher{Springer}, \baddress{New York}.
\bid{doi={10.1007/978-1-4419-0661-8}, mr={2567175}}
\bptok{imsref}%
\end{bbook}
%
\endbibitem

\bibitem[\protect\citeauthoryear{Bai and Yin}{1993}]{bai1993}
%
\begin{barticle}[mr]
\bauthor{\bsnm{Bai},~\bfnm{Z.~D.}\binits{Z.~D.}} \AND
\bauthor{\bsnm{Yin},~\bfnm{Y.~Q.}\binits{Y.~Q.}}
(\byear{1993}).
\btitle{Limit of the smallest eigenvalue of a~large-dimensional sample
covariance matrix}.
\bjournal{Ann. Probab.}
\bvolume{21}
\bpages{1275--1294}.
\bid{issn={0091-1798}, mr={1235416}}
\bptok{imsref}%
\end{barticle}
%
\endbibitem

\bibitem[\protect\citeauthoryear{Bentkus and Rudzkis}{1976}]{bentkus1976}
%
\begin{barticle}[mr]
\bauthor{\bsnm{Bentkus},~\bfnm{R.}\binits{R.}} \AND
\bauthor{\bsnm{Rudzkis},~\bfnm{R.}\binits{R.}}
(\byear{1976}).
\btitle{Large deviations for estimates of the spectrum of a~stationary
{G}aussian sequence}.
\bjournal{Litovsk. Mat. Sb.}
\bvolume{16}
\bpages{63--77, 253}.
\bid{issn={0132-2818}, mr={0436510}}
\bptok{imsref}%
\end{barticle}
%
\endbibitem

\bibitem[\protect\citeauthoryear{Bercu, Gamboa and Rouault}{1997}]{bercu1997}
%
\begin{barticle}[mr]
\bauthor{\bsnm{Bercu},~\bfnm{B.}\binits{B.}},
\bauthor{\bsnm{Gamboa},~\bfnm{F.}\binits{F.}} \AND
\bauthor{\bsnm{Rouault},~\bfnm{A.}\binits{A.}}
(\byear{1997}).
\btitle{Large deviations for quadratic forms of stationary {G}aussian
processes}.
\bjournal{Stochastic Process. Appl.}
\bvolume{71}
\bpages{75--90}.
\bid{doi={10.1016/S0304-4149(97)00071-9}, issn={0304-4149}, mr={1480640}}
\bptok{imsref}%
\end{barticle}
%
\endbibitem

\bibitem[\protect\citeauthoryear{Bercu, Gamboa and
Lavielle}{2000}]{bercu2000}
%
\begin{barticle}[mr]
\bauthor{\bsnm{Bercu},~\bfnm{Bernard}\binits{B.}},
\bauthor{\bsnm{Gamboa},~\bfnm{Fabrice}\binits{F.}} \AND
\bauthor{\bsnm{Lavielle},~\bfnm{Marc}\binits{M.}}
(\byear{2000}).
\btitle{Sharp large deviations for {G}aussian quadratic forms with
applications}.
\bjournal{ESAIM Probab. Stat.}
\bvolume{4}
\bpages{1--24 (electronic)}.
\bid{doi={10.1051/ps:2000101}, issn={1292-8100}, mr={1749403}}
\bptok{imsref}%
\end{barticle}
%
\endbibitem

\bibitem[\protect\citeauthoryear{Bickel and Levina}{2008a}]{bickel2008b}
%
\begin{barticle}[mr]
\bauthor{\bsnm{Bickel},~\bfnm{Peter~J.}\binits{P.~J.}} \AND
\bauthor{\bsnm{Levina},~\bfnm{Elizaveta}\binits{E.}}
(\byear{2008}a).
\btitle{Covariance regularization by thresholding}.
\bjournal{Ann. Statist.}
\bvolume{36}
\bpages{2577--2604}.
\bid{doi={10.1214/08-AOS600}, issn={0090-5364}, mr={2485008}}
\bptok{imsref}%
\end{barticle}
%
\endbibitem

\bibitem[\protect\citeauthoryear{Bickel and Levina}{2008b}]{bickel2008a}
%
\begin{barticle}[mr]
\bauthor{\bsnm{Bickel},~\bfnm{Peter~J.}\binits{P.~J.}} \AND
\bauthor{\bsnm{Levina},~\bfnm{Elizaveta}\binits{E.}}
(\byear{2008}b).
\btitle{Regularized estimation of large covariance matrices}.
\bjournal{Ann. Statist.}
\bvolume{36}
\bpages{199--227}.
\bid{doi={10.1214/009053607000000758}, issn={0090-5364}, mr={2387969}}
\bptok{imsref}%
\end{barticle}
%
\endbibitem

\bibitem[\protect\citeauthoryear{Bryc and Dembo}{1997}]{bryc1997}
%
\begin{barticle}[mr]
\bauthor{\bsnm{Bryc},~\bfnm{W{\l}odzimierz}\binits{W.}} \AND
\bauthor{\bsnm{Dembo},~\bfnm{Amir}\binits{A.}}
(\byear{1997}).
\btitle{Large deviations for quadratic functionals of {G}aussian processes}.
\bjournal{J. Theoret. Probab.}
\bvolume{10}
\bpages{307--332}.
\bnote{Dedicated to Murray Rosenblatt}.
\bid{doi={10.1023/A:1022656331883}, issn={0894-9840}, mr={1455147}}
\bptok{imsref}%
\end{barticle}
%
\endbibitem

\bibitem[\protect\citeauthoryear{Bryc, Dembo and Jiang}{2006}]{bryc2006}
%
\begin{barticle}[mr]
\bauthor{\bsnm{Bryc},~\bfnm{W{\l}odzimierz}\binits{W.}},
\bauthor{\bsnm{Dembo},~\bfnm{Amir}\binits{A.}} \AND
\bauthor{\bsnm{Jiang},~\bfnm{Tiefeng}\binits{T.}}
(\byear{2006}).
\btitle{Spectral measure of large random {H}ankel, {M}arkov and {T}oeplitz
matrices}.
\bjournal{Ann. Probab.}
\bvolume{34}
\bpages{1--38}.
\bid{doi={10.1214/009117905000000495}, issn={0091-1798}, mr={2206341}}
\bptok{imsref}%
\end{barticle}
%
\endbibitem

\bibitem[\protect\citeauthoryear{B{\"u}hlmann and
K{\"u}nsch}{1999}]{buhlmann1999}
%
\begin{barticle}[author]
\bauthor{\bsnm{B{\"u}hlmann},~\bfnm{P.}\binits{P.}} \AND
\bauthor{\bsnm{K{\"u}nsch},~\bfnm{H.~R.}\binits{H.~R.}}
(\byear{1999}).
\btitle{Block length selection in the bootstrap for time series}.
\bjournal{Comput. Statist. Data Anal.}
\bvolume{31}
\bpages{295--310}.
\bptok{imsref}%
\end{barticle}
%
\endbibitem

\bibitem[\protect\citeauthoryear{Burkholder}{1988}]{burkholder1988}
%
\begin{barticle}[mr]
\bauthor{\bsnm{Burkholder},~\bfnm{Donald~L.}\binits{D.~L.}}
(\byear{1988}).
\btitle{Sharp inequalities for martingales and stochastic integrals}.
\bjournal{Colloque Paul L\'{e}vy sur les Processus Stochastiques
(Palaiseau, 1987).
Ast\'erisque}
\bvolume{157-158}
\bpages{75--94}.
\bid{issn={0303-1179}, mr={0976214}}
\bptok{imsref}%
\end{barticle}
%
\endbibitem

\bibitem[\protect\citeauthoryear{Cai, Zhang and Zhou}{2010}]{cai2010}
%
\begin{barticle}[mr]
\bauthor{\bsnm{Cai},~\bfnm{T.~Tony}\binits{T.~T.}},
\bauthor{\bsnm{Zhang},~\bfnm{Cun-Hui}\binits{C.-H.}} \AND
\bauthor{\bsnm{Zhou},~\bfnm{Harrison~H.}\binits{H.~H.}}
(\byear{2010}).
\btitle{Optimal rates of convergence for covariance matrix estimation}.
\bjournal{Ann. Statist.}
\bvolume{38}
\bpages{2118--2144}.
\bid{doi={10.1214/09-AOS752}, issn={0090-5364}, mr={2676885}}
\bptok{imsref}%
\end{barticle}
%
\endbibitem

\bibitem[\protect\citeauthoryear{Dembo and Zeitouni}{1998}]{dembo1998}
%
\begin{bbook}[mr]
\bauthor{\bsnm{Dembo},~\bfnm{Amir}\binits{A.}} \AND
\bauthor{\bsnm{Zeitouni},~\bfnm{Ofer}\binits{O.}}
(\byear{1998}).
\btitle{Large Deviations Techniques and Applications},
\bedition{2nd} ed.
\bseries{Applications of Mathematics (New York)}
\bvolume{38}.
\bpublisher{Springer}, \baddress{New York}.
\bid{mr={1619036}}
\bptok{imsref}%
\end{bbook}
%
\endbibitem

\bibitem[\protect\citeauthoryear{Djellout, Guillin and
Wu}{2006}]{djellout2006}
%
\begin{barticle}[mr]
\bauthor{\bsnm{Djellout},~\bfnm{H.}\binits{H.}},
\bauthor{\bsnm{Guillin},~\bfnm{A.}\binits{A.}} \AND
\bauthor{\bsnm{Wu},~\bfnm{L.}\binits{L.}}
(\byear{2006}).
\btitle{Moderate deviations of empirical periodogram and non-linear functionals
of moving average processes}.
\bjournal{Ann. Inst. Henri Poincar\'e Probab. Stat.}
\bvolume{42}
\bpages{393--416}.
\bid{doi={10.1016/j.anihpb.2005.04.006}, issn={0246-0203}, mr={2242954}}
\bptok{imsref}%
\end{barticle}
%
\endbibitem

\bibitem[\protect\citeauthoryear{El~Karoui}{2005}]{karoui2005}
%
\begin{barticle}[mr]
\bauthor{\bsnm{El~Karoui},~\bfnm{Noureddine}\binits{N.}}
(\byear{2005}).
\btitle{Recent results about the largest eigenvalue of random covariance
matrices and statistical application}.
\bjournal{Acta Phys. Polon. B}
\bvolume{36}
\bpages{2681--2697}.
\bid{issn={0587-4254}, mr={2188088}}
\bptok{imsref}%
\end{barticle}
%
\endbibitem

\bibitem[\protect\citeauthoryear{El~Karoui}{2008}]{karoui2008}
%
\begin{barticle}[mr]
\bauthor{\bsnm{El~Karoui},~\bfnm{Noureddine}\binits{N.}}
(\byear{2008}).
\btitle{Operator norm consistent estimation of large-dimensional sparse
covariance matrices}.
\bjournal{Ann. Statist.}
\bvolume{36}
\bpages{2717--2756}.
\bid{doi={10.1214/07-AOS559}, issn={0090-5364}, mr={2485011}}
\bptok{imsref}%
\end{barticle}
%
\endbibitem

\bibitem[\protect\citeauthoryear{Freedman}{1975}]{freedman1975}
%
\begin{barticle}[mr]
\bauthor{\bsnm{Freedman},~\bfnm{David~A.}\binits{D.~A.}}
(\byear{1975}).
\btitle{On tail probabilities for martingales}.
\bjournal{Ann. Probab.}
\bvolume{3}
\bpages{100--118}.
\bid{mr={0380971}}
\bptok{imsref}%
\end{barticle}
%
\endbibitem

\bibitem[\protect\citeauthoryear{Furrer and Bengtsson}{2007}]{furrer2007}
%
\begin{barticle}[mr]
\bauthor{\bsnm{Furrer},~\bfnm{Reinhard}\binits{R.}} \AND
\bauthor{\bsnm{Bengtsson},~\bfnm{Thomas}\binits{T.}}
(\byear{2007}).
\btitle{Estimation of high-dimensional prior and posterior covariance matrices
in {K}alman filter variants}.
\bjournal{J. Multivariate Anal.}
\bvolume{98}
\bpages{227--255}.
\bid{doi={10.1016/j.jmva.2006.08.003}, issn={0047-259X}, mr={2301751}}
\bptok{imsref}%
\end{barticle}
%
\endbibitem

\bibitem[\protect\citeauthoryear{Gamboa, Rouault and Zani}{1999}]{gamboa1999}
%
\begin{barticle}[mr]
\bauthor{\bsnm{Gamboa},~\bfnm{F.}\binits{F.}},
\bauthor{\bsnm{Rouault},~\bfnm{A.}\binits{A.}} \AND
\bauthor{\bsnm{Zani},~\bfnm{M.}\binits{M.}}
(\byear{1999}).
\btitle{A~functional large deviations principle for quadratic forms of
{G}aussian stationary processes}.
\bjournal{Statist. Probab. Lett.}
\bvolume{43}
\bpages{299--308}.
\bid{doi={10.1016/S0167-7152(98)00270-3}, issn={0167-7152}, mr={1708097}}
\bptok{imsref}%
\end{barticle}
%
\endbibitem

\bibitem[\protect\citeauthoryear{Geman}{1980}]{geman1980}
%
\begin{barticle}[mr]
\bauthor{\bsnm{Geman},~\bfnm{Stuart}\binits{S.}}
(\byear{1980}).
\btitle{A~limit theorem for the norm of random matrices}.
\bjournal{Ann. Probab.}
\bvolume{8}
\bpages{252--261}.
\bid{issn={0091-1798}, mr={0566592}}
\bptok{imsref}%
\end{barticle}
%
\endbibitem

\bibitem[\protect\citeauthoryear{Grenander and Szeg{\"
o}}{1958}]{grenander1958}
%
\begin{bbook}[mr]
\bauthor{\bsnm{Grenander},~\bfnm{Ulf}\binits{U.}} \AND
\bauthor{\bsnm{Szeg{\"o}},~\bfnm{Gabor}\binits{G.}}
(\byear{1958}).
\btitle{Toeplitz Forms and Their Applications}.
\bpublisher{Univ. California Press}, \baddress{Berkeley}.
\bid{mr={0094840}}
\bptok{imsref}%
\end{bbook}
%
\endbibitem

\bibitem[\protect\citeauthoryear{Haeusler}{1984}]{haeusler1984}
%
\begin{barticle}[mr]
\bauthor{\bsnm{Haeusler},~\bfnm{Erich}\binits{E.}}
(\byear{1984}).
\btitle{An exact rate of convergence in the functional central limit theorem
for special martingale difference arrays}.
\bjournal{Z. Wahrsch. Verw. Gebiete}
\bvolume{65}
\bpages{523--534}.
\bid{doi={10.1007/BF00531837}, issn={0044-3719}, mr={0736144}}
\bptok{imsref}%
\end{barticle}
%
\endbibitem

\bibitem[\protect\citeauthoryear{Horn and Johnson}{1990}]{horn1990}
%
\begin{bbook}[mr]
\bauthor{\bsnm{Horn},~\bfnm{Roger~A.}\binits{R.~A.}} \AND
\bauthor{\bsnm{Johnson},~\bfnm{Charles~R.}\binits{C.~R.}}
(\byear{1990}).
\btitle{Matrix Analysis}.
\bpublisher{Cambridge Univ. Press}, \baddress{Cambridge}.
\bnote{Corrected reprint of the 1985 original}.
\bid{mr={1084815}}
\bptok{imsref}%
\end{bbook}
%
\endbibitem

\bibitem[\protect\citeauthoryear{Johansson}{2000}]{johansson2000}
%
\begin{barticle}[mr]
\bauthor{\bsnm{Johansson},~\bfnm{Kurt}\binits{K.}}
(\byear{2000}).
\btitle{Shape fluctuations and random matrices}.
\bjournal{Comm. Math. Phys.}
\bvolume{209}
\bpages{437--476}.
\bid{doi={10.1007/s002200050027}, issn={0010-3616}, mr={1737991}}
\bptok{imsref}%
\end{barticle}
%
\endbibitem

\bibitem[\protect\citeauthoryear{Johnstone}{2001}]{johnstone2001}
%
\begin{barticle}[mr]
\bauthor{\bsnm{Johnstone},~\bfnm{Iain~M.}\binits{I.~M.}}
(\byear{2001}).
\btitle{On the distribution of the largest eigenvalue in principal components
analysis}.
\bjournal{Ann. Statist.}
\bvolume{29}
\bpages{295--327}.
\bid{doi={10.1214/aos/1009210544}, issn={0090-5364}, mr={1863961}}
\bptok{imsref}%
\end{barticle}
%
\endbibitem

\bibitem[\protect\citeauthoryear{Kakizawa}{2007}]{kakizawa2007}
%
\begin{barticle}[mr]
\bauthor{\bsnm{Kakizawa},~\bfnm{Yoshihide}\binits{Y.}}
(\byear{2007}).
\btitle{Moderate deviations for quadratic forms in {G}aussian stationary
processes}.
\bjournal{J.~Multivariate Anal.}
\bvolume{98}
\bpages{992--1017}.
\bid{doi={10.1016/j.jmva.2006.07.004}, issn={0047-259X}, mr={2325456}}
\bptok{imsref}%
\end{barticle}
%
\endbibitem

\bibitem[\protect\citeauthoryear{Kolmogoroff}{1941}]{kolmogorov1941}
%
\begin{barticle}[mr]
\bauthor{\bsnm{Kolmogoroff},~\bfnm{A.}\binits{A.}}
(\byear{1941}).
\btitle{Interpolation und {E}xtrapolation von station\"aren zuf\"alligen
{F}olgen}.
\bjournal{Bull. Acad. Sci. URSS S\'er. Math. [Izvestia Akad. Nauk SSSR]}
\bvolume{5}
\bpages{3--14}.
\bid{mr={0004416}}
\bptok{imsref}%
\end{barticle}
%
\endbibitem

\bibitem[\protect\citeauthoryear{Lin and Liu}{2009}]{lin2009}
%
\begin{barticle}[mr]
\bauthor{\bsnm{Lin},~\bfnm{Zhengyan}\binits{Z.}} \AND
\bauthor{\bsnm{Liu},~\bfnm{Weidong}\binits{W.}}
(\byear{2009}).
\btitle{On maxima of periodograms of stationary processes}.
\bjournal{Ann. Statist.}
\bvolume{37}
\bpages{2676--2695}.
\bid{doi={10.1214/08-AOS590}, issn={0090-5364}, mr={2541443}}
\bptok{imsref}%
\end{barticle}
%
\endbibitem

\bibitem[\protect\citeauthoryear{Liu and Shao}{2010}]{liu2010}
%
\begin{barticle}[mr]
\bauthor{\bsnm{Liu},~\bfnm{Weidong}\binits{W.}} \AND
\bauthor{\bsnm{Shao},~\bfnm{Qi-Man}\binits{Q.-M.}}
(\byear{2010}).
\btitle{Cram\'er-type moderate deviation for the maximum of the periodogram
with application to simultaneous tests in gene expression time series}.
\bjournal{Ann. Statist.}
\bvolume{38}
\bpages{1913--1935}.
\bid{doi={10.1214/09-AOS774}, issn={0090-5364}, mr={2662363}}
\bptok{imsref}%
\end{barticle}
%
\endbibitem

\bibitem[\protect\citeauthoryear{Liu and Wu}{2010}]{liuwu2010}
%
\begin{barticle}[mr]
\bauthor{\bsnm{Liu},~\bfnm{Weidong}\binits{W.}} \AND
\bauthor{\bsnm{Wu},~\bfnm{Wei~Biao}\binits{W.~B.}}
(\byear{2010}).
\btitle{Asymptotics of spectral density estimates}.
\bjournal{Econometric Theory}
\bvolume{26}
\bpages{1218--1245}.
\bid{doi={10.1017/S026646660999051X}, issn={0266-4666}, mr={2660298}}
\bptok{imsref}%
\end{barticle}
%
\endbibitem

\bibitem[\protect\citeauthoryear{Mar{\v{c}}enko and
Pastur}{1967}]{marcenko1967}
%
\begin{barticle}[mr]
\bauthor{\bsnm{Mar{\v{c}}enko},~\bfnm{V.~A.}\binits{V.~A.}} \AND
\bauthor{\bsnm{Pastur},~\bfnm{L.~A.}\binits{L.~A.}}
(\byear{1967}).
\btitle{Distribution of eigenvalues in certain sets of random matrices}.
\bjournal{Mat. Sb. (N.S.)}
\bvolume{72}
\bpages{507--536}.
\bid{mr={0208649}}
\bptok{imsref}%
\end{barticle}
%
\endbibitem

\bibitem[\protect\citeauthoryear{McMurry and Politis}{2010}]{mcmurry2010}
%
\begin{barticle}[mr]
\bauthor{\bsnm{McMurry},~\bfnm{Timothy~L.}\binits{T.~L.}} \AND
\bauthor{\bsnm{Politis},~\bfnm{Dimitris~N.}\binits{D.~N.}}
(\byear{2010}).
\btitle{Banded and tapered estimates for autocovariance matrices and
the linear
process bootstrap}.
\bjournal{J. Time Series Anal.}
\bvolume{31}
\bpages{471--482}.
\bid{doi={10.1111/j.1467-9892.2010.00679.x}, issn={0143-9782}, mr={2732601}}
\bptok{imsref}%
\end{barticle}
%
\endbibitem

\bibitem[\protect\citeauthoryear{Nagaev}{1979}]{nagaev1979}
%
\begin{barticle}[mr]
\bauthor{\bsnm{Nagaev},~\bfnm{S.~V.}\binits{S.~V.}}
(\byear{1979}).
\btitle{Large deviations of sums of independent random variables}.
\bjournal{Ann. Probab.}
\bvolume{7}
\bpages{745--789}.
\bid{issn={0091-1798}, mr={0542129}}
\bptok{imsref}%
\end{barticle}
%
\endbibitem

\bibitem[\protect\citeauthoryear{Peligrad and Wu}{2010}]{peligrad2010}
%
\begin{barticle}[mr]
\bauthor{\bsnm{Peligrad},~\bfnm{Magda}\binits{M.}} \AND
\bauthor{\bsnm{Wu},~\bfnm{Wei~Biao}\binits{W.~B.}}
(\byear{2010}).
\btitle{Central limit theorem for {F}ourier transforms of stationary
processes}.
\bjournal{Ann. Probab.}
\bvolume{38}
\bpages{2009--2022}.
\bid{doi={10.1214/10-AOP530}, issn={0091-1798}, mr={2722793}}
\bptok{imsref}%
\end{barticle}
%
\endbibitem

\bibitem[\protect\citeauthoryear{Politis}{2003}]{politis2003}
%
\begin{barticle}[mr]
\bauthor{\bsnm{Politis},~\bfnm{Dimitris~N.}\binits{D.~N.}}
(\byear{2003}).
\btitle{Adaptive bandwidth choice}.
\bjournal{J. Nonparametr. Stat.}
\bvolume{15}
\bpages{517--533}.
\bid{doi={10.1080/10485250310001604659}, issn={1048-5252}, mr={2017485}}
\bptok{imsref}%
\end{barticle}
%
\endbibitem

\bibitem[\protect\citeauthoryear{Politis, Romano and
Wolf}{1999}]{politis1999}
%
\begin{bbook}[mr]
\bauthor{\bsnm{Politis},~\bfnm{Dimitris~N.}\binits{D.~N.}},
\bauthor{\bsnm{Romano},~\bfnm{Joseph~P.}\binits{J.~P.}} \AND
\bauthor{\bsnm{Wolf},~\bfnm{Michael}\binits{M.}}
(\byear{1999}).
\btitle{Subsampling}.
\bpublisher{Springer}, \baddress{New York}.
\bid{mr={1707286}}
\bptok{imsref}%
\end{bbook}
%
\endbibitem

\bibitem[\protect\citeauthoryear{Rio}{2009}]{rio2009}
%
\begin{barticle}[mr]
\bauthor{\bsnm{Rio},~\bfnm{Emmanuel}\binits{E.}}
(\byear{2009}).
\btitle{Moment inequalities for sums of dependent random variables under
projective conditions}.
\bjournal{J. Theoret. Probab.}
\bvolume{22}
\bpages{146--163}.
\bid{doi={10.1007/s10959-008-0155-9}, issn={0894-9840}, mr={2472010}}
\bptok{imsref}%
\end{barticle}
%
\endbibitem

\bibitem[\protect\citeauthoryear{Rudzkis}{1978}]{rudzkis1978}
%
\begin{barticle}[mr]
\bauthor{\bsnm{Rudzkis},~\bfnm{R.}\binits{R.}}
(\byear{1978}).
\btitle{Large deviations for estimates of the spectrum of a~stationary
sequence}.
\bjournal{Litovsk. Mat. Sb.}
\bvolume{18}
\bpages{81--98, 217}.
\bid{issn={0132-2818}, mr={0519099}}
\bptok{imsref}%
\end{barticle}
%
\endbibitem

\bibitem[\protect\citeauthoryear{Saulis and
Statulevi{\v{c}}ius}{1991}]{saulis1991}
%
\begin{bbook}[mr]
\bauthor{\bsnm{Saulis},~\bfnm{L.}\binits{L.}} \AND
\bauthor{\bsnm{Statulevi{\v{c}}ius},~\bfnm{V.~A.}\binits{V.~A.}}
(\byear{1991}).
\btitle{Limit Theorems for Large Deviations}.
\bseries{Mathematics and Its Applications (Soviet Series)}
\bvolume{73}.
\bpublisher{Kluwer}, \baddress{Dordrecht}.
\bnote{Translated and revised from the 1989 Russian original}.
\bid{mr={1171883}}
\bptok{imsref}%
\end{bbook}
%
\endbibitem

\bibitem[\protect\citeauthoryear{Shao and Wu}{2007}]{shao2007}
%
\begin{barticle}[mr]
\bauthor{\bsnm{Shao},~\bfnm{Xiaofeng}\binits{X.}} \AND
\bauthor{\bsnm{Wu},~\bfnm{Wei~Biao}\binits{W.~B.}}
(\byear{2007}).
\btitle{Asymptotic spectral theory for nonlinear time series}.
\bjournal{Ann. Statist.}
\bvolume{35}
\bpages{1773--1801}.
\bid{doi={10.1214/009053606000001479}, issn={0090-5364}, mr={2351105}}
\bptok{imsref}%
\end{barticle}
%
\endbibitem

\bibitem[\protect\citeauthoryear{Solo}{2010}]{solo2010}
%
\begin{bmisc}[author]
\bauthor{\bsnm{Solo},~\bfnm{V.}\binits{V.}}
(\byear{2010}).
\bhowpublished{On random matrix theory for stationary processes.
In \textit{IEEE International Conference on}
\textit{Acoustics Speech and Signal Processing} (\textit{ICASSP})
3758--3761. IEEE, Piscataway, NJ.}
\bptok{imsref}%
\end{bmisc}
%
\endbibitem

\bibitem[\protect\citeauthoryear{Toeplitz}{1911}]{toeplitz1911}
%
\begin{barticle}[mr]
\bauthor{\bsnm{Toeplitz},~\bfnm{Otto}\binits{O.}}
(\byear{1911}).
\btitle{Zur {T}heorie der quadratischen und bilinearen {F}ormen von
unendlichvielen {V}er\"anderlichen}.
\bjournal{Math. Ann.}
\bvolume{70}
\bpages{351--376}.
\bid{doi={10.1007/BF01564502}, issn={0025-5831}, mr={1511625}}
\bptok{imsref}%
\end{barticle}
%
\endbibitem

\bibitem[\protect\citeauthoryear{Tong}{1990}]{tong1990}
%
\begin{bbook}[mr]
\bauthor{\bsnm{Tong},~\bfnm{Howell}\binits{H.}}
(\byear{1990}).
\btitle{Nonlinear Time Series}.
\bseries{Oxford Statistical Science Series}
\bvolume{6}.
\bpublisher{Oxford Univ. Press}, \baddress{New York}.
\bid{mr={1079320}}
\bptok{imsref}%
\end{bbook}
%
\endbibitem

\bibitem[\protect\citeauthoryear{Tracy and Widom}{1994}]{tracy1994}
%
\begin{barticle}[mr]
\bauthor{\bsnm{Tracy},~\bfnm{Craig~A.}\binits{C.~A.}} \AND
\bauthor{\bsnm{Widom},~\bfnm{Harold}\binits{H.}}
(\byear{1994}).
\btitle{Level-spacing distributions and the {A}iry kernel}.
\bjournal{Comm. Math. Phys.}
\bvolume{159}
\bpages{151--174}.
\bid{issn={0010-3616}, mr={1257246}}
\bptok{imsref}%
\end{barticle}
%
\endbibitem

\bibitem[\protect\citeauthoryear{Turkman and Walker}{1984}]{turkman1984}
%
\begin{barticle}[mr]
\bauthor{\bsnm{Turkman},~\bfnm{K.~F.}\binits{K.~F.}} \AND
\bauthor{\bsnm{Walker},~\bfnm{A.~M.}\binits{A.~M.}}
(\byear{1984}).
\btitle{On the asymptotic distributions of maxima of trigonometric polynomials
with random coefficients}.
\bjournal{Adv. in Appl. Probab.}
\bvolume{16}
\bpages{819--842}.
\bid{doi={10.2307/1427342}, issn={0001-8678}, mr={0766781}}
\bptok{imsref}%
\end{barticle}
%
\endbibitem

\bibitem[\protect\citeauthoryear{Turkman and Walker}{1990}]{turkman1990}
%
\begin{barticle}[mr]
\bauthor{\bsnm{Turkman},~\bfnm{K.~F.}\binits{K.~F.}} \AND
\bauthor{\bsnm{Walker},~\bfnm{A.~M.}\binits{A.~M.}}
(\byear{1990}).
\btitle{A~stability result for the periodogram}.
\bjournal{Ann. Probab.}
\bvolume{18}
\bpages{1765--1783}.
\bid{issn={0091-1798}, mr={1071824}}
\bptok{imsref}%
\end{barticle}
%
\endbibitem

\bibitem[\protect\citeauthoryear{Wiener}{1949}]{wiener1949}
%
\begin{bbook}[mr]
\bauthor{\bsnm{Wiener},~\bfnm{Norbert}\binits{N.}}
(\byear{1949}).
\btitle{Extrapolation, {I}nterpolation, and {S}moothing of
Stationary Time Series. {W}ith {E}ngineering {A}pplications}.
\bpublisher{MIT Press},
\baddress{Cambridge, MA}.
\bid{mr={0031213}}
\bptok{imsref}%
\end{bbook}
%
\endbibitem

\bibitem[\protect\citeauthoryear{Woodroofe and
Van~Ness}{1967}]{woodroofe1967}
%
\begin{barticle}[mr]
\bauthor{\bsnm{Woodroofe},~\bfnm{Michael~B.}\binits{M.~B.}} \AND
\bauthor{\bsnm{Van~Ness},~\bfnm{John~W.}\binits{J.~W.}}
(\byear{1967}).
\btitle{The maximum deviation of sample spectral densities}.
\bjournal{Ann. Math. Statist.}
\bvolume{38}
\bpages{1558--1569}.
\bid{issn={0003-4851}, mr={0216717}}
\bptok{imsref}%
\end{barticle}
%
\endbibitem

\bibitem[\protect\citeauthoryear{Wu}{2005}]{wu2005}
%
\begin{barticle}[mr]
\bauthor{\bsnm{Wu},~\bfnm{Wei~Biao}\binits{W.~B.}}
(\byear{2005}).
\btitle{Nonlinear system theory: Another look at dependence}.
\bjournal{Proc. Natl. Acad. Sci. USA}
\bvolume{102}
\bpages{14150--14154 (electronic)}.
\bid{doi={10.1073/pnas.0506715102}, issn={1091-6490}, mr={2172215}}
\bptok{imsref}%
\end{barticle}
%
\endbibitem

\bibitem[\protect\citeauthoryear{Wu and Pourahmadi}{2009}]{wuss2009}
%
\begin{barticle}[mr]
\bauthor{\bsnm{Wu},~\bfnm{Wei~Biao}\binits{W.~B.}} \AND
\bauthor{\bsnm{Pourahmadi},~\bfnm{Mohsen}\binits{M.}}
(\byear{2009}).
\btitle{Banding sample autocovariance matrices of stationary processes}.
\bjournal{Statist. Sinica}
\bvolume{19}
\bpages{1755--1768}.
\bid{issn={1017-0405}, mr={2589209}}
\bptok{imsref}%
\end{barticle}
%
\endbibitem

\bibitem[\protect\citeauthoryear{Xiao and Wu}{2011}]{wu2011}
%
\begin{bmisc}[author]
\bauthor{\bsnm{Xiao},~\bfnm{H.}\binits{H.}} \AND
\bauthor{\bsnm{Wu},~\bfnm{W.~B.}\binits{W.~B.}}
(\byear{2011}).
\bhowpublished{Asymptotic inference of autocovariances of stationary processes.
Available at arXiv:\arxivurl{1105.3423}.}
\bptok{imsref}%
\end{bmisc}
%
\endbibitem

\bibitem[\protect\citeauthoryear{Xiao and Wu}{2012}]{wu2011s}
%
\begin{bmisc}[author]
\bauthor{\bsnm{Xiao},~\bfnm{H}\binits{H.}} \AND
\bauthor{\bsnm{Wu},~\bfnm{W.~B.}\binits{W.~B.}}
(\byear{2012}).
\bhowpublished{Supplement to ``Covariance matrix estimation
for stationary time series.''
DOI:\doiurl{10.1214/11-AOS967SUPP}.}
\bptok{imsref}%
\end{bmisc}
%
\endbibitem

\bibitem[\protect\citeauthoryear{Yin, Bai and Krishnaiah}{1988}]{yin1988}
%
\begin{barticle}[mr]
\bauthor{\bsnm{Yin},~\bfnm{Y.~Q.}\binits{Y.~Q.}},
\bauthor{\bsnm{Bai},~\bfnm{Z.~D.}\binits{Z.~D.}} \AND
\bauthor{\bsnm{Krishnaiah},~\bfnm{P.~R.}\binits{P.~R.}}
(\byear{1988}).
\btitle{On the limit of the largest eigenvalue of the
large-dimensional sample
covariance matrix}.
\bjournal{Probab. Theory Related Fields}
\bvolume{78}
\bpages{509--521}.
\bid{doi={10.1007/BF00353874}, issn={0178-8051}, mr={0950344}}
\bptok{imsref}%
\end{barticle}
%
\endbibitem

\bibitem[\protect\citeauthoryear{Zani}{2002}]{zani2002}
%
\begin{barticle}[mr]
\bauthor{\bsnm{Zani},~\bfnm{Marguerite}\binits{M.}}
(\byear{2002}).
\btitle{Large deviations for quadratic forms of locally stationary processes}.
\bjournal{J.~Multivariate Anal.}
\bvolume{81}
\bpages{205--228}.
\bid{doi={10.1006/jmva.2001.2003}, issn={0047-259X}, mr={1906377}}
\bptok{imsref}%
\end{barticle}
%
\endbibitem

\bibitem[\protect\citeauthoryear{Zygmund}{2002}]{zygmund2002}
%
\begin{bbook}[mr]
\bauthor{\bsnm{Zygmund},~\bfnm{A.}\binits{A.}}
(\byear{2002}).
\btitle{Trigonometric Series. {V}ols {I}, {II}},
\bedition{3rd} ed.
\bpublisher{Cambridge Univ. Press}, \baddress{Cambridge}.
\bid{mr={1963498}}
\bptok{imsref}%
\end{bbook}
%
\endbibitem

\end{thebibliography}
\end{document}